\documentclass[12pt]{article}
\def\newpic#1{%
\def\emline##1##2##3##4##5##6{%
\put(##1,##2){\special{em:point #1##3}}%
\put(##4,##5){\special{em:point #1##6}}%
\special{em:line #1##3,#1##6}}}
\newpic{}
\def\emline#1#2#3#4#5#6{%
\put(#1,#2){\special{em:moveto}}%
\put(#4,#5){\special{em:lineto}}}
\def\newpic#1{}
\newcommand\NN{\hbox{I\kern-.2em\hbox{N}}}
\newcommand\RR{\hbox{I\kern-.2em\hbox{R}}}
\newcommand\sRR{{\sl \hbox{I\kern-.2em\hbox{R}}}}
\newcommand\QQ{\hbox{I\kern-.53em\hbox{Q}}}
\newcommand\ZZ{{{\rm Z}\kern-.28em{\rm Z}}}

\title{On a $\{K_4,K_{2,2,2}\}$-ultrahomogeneous graph}
\author{Italo J. Dejter
\\ University of Puerto Rico \\ Rio Piedras, PR 00931-3355 \\ ijdejter@uprrp.edu}
\newtheorem{lem}{Lemma}[section]
\newtheorem{teo}[lem]{Theorem}

\newenvironment{demo}{\noindent {\bf Proof:}\ }{\ \framebox{}}

\newcommand{\example}{\noindent{\bf Example.\hspace{4mm}}}

\date{}
\begin{document}
\maketitle

\begin{abstract}
The existence of a connected 12-regular
$\{K_4,K_{2,2,2}\}$-ul\-tra\-ho\-mo\-ge\-neous graph $G$ is
established, (i.e. each isomorphism between two copies of $K_4$ or
$K_{2,2,2}$ in $G$ extends to an automorphism of $G$), with the 42
ordered lines of the Fano plane taken as vertices. This graph $G$ can be
expressed in a unique way both as the edge-disjoint union of 42
induced copies of $K_4$ and as the edge-disjoint union of 21 induced
copies of $K_{2,2,2}$, with no more copies of $K_4$ or $K_{2,2,2}$
existing in $G$. Moreover, each edge of $G$ is shared by exactly one
copy of $K_4$ and one of $K_{2,2,2}$. While the line graphs of
$d$-cubes, ($3\le d\in\ZZ$), are $\{K_d,
K_{2,2}\}$-ultrahomogeneous, $G$ is not even line-graphical. In
addition, the chordless 6-cycles of $G$ are seen to play an interesting role and some self-dual configurations associated to $G$ with 2-arc-transitive, arc-transitive and semisymmetric Levi graphs are considered.
\end{abstract}

\section{Introduction}

Let $H$ be a connected regular graph and let $m,n\in\ZZ$ with
$1<m<n$. An $\{H\}_n^m$-{\it graph} is a connected graph
that: {\bf(a)} is representable as an edge-disjoint union
of $n$ induced copies of $H$; {\bf(b)} has exactly $m$ copies of $H$
incident to each vertex, with no two such copies sharing more than
one vertex; and {\bf(c)} has exactly $n$ copies of $H$ as
induced subgraphs isomorphic to $H$.

We remark that an $\{H\}_n^m$-graph $G$ is $\{H\}$-{\it ultrahomogeneous} (as in
\cite{I}) if every isomorphism between two copies of $H$ in $G$
extends to an auto\-mor\-phism of $G$. Graph ultrahomogeneity is a
concept that can be traced back to \cite{Sheh,Gard,Ronse}.

Notice that a connected graph $G$ is $m$-regular if and
only if it is a $\{K_2\}_{|E(G)|}^m$-graph. In this case, $G$ is
arc-transitive if and only if $G$ is $\{K_2\}$-ultrahomogeneous.
Thus, $\{H\}$-ultrahomogeneity is a notion of graph symmetry stronger than
arc-transitivity.

If $G$ is an $\{H_i\}_{n_i}^{m_i}$-graph, where $i=1,2$, and $H_1\ne
H_2$, then $G$ is said to be an $\{H_1\}_{n_1}^{m_1}\{H_2\}_{n_2}^{m_2}$-{\it graph}.
If, in  addition,
$G$ is $\{H_i\}$-ultrahomogeneous, for both $i=1,2$, then $G$
is $\{H_1,H_2\}$-{\it ultrahomogeneous}, again as in \cite{I}.
If each edge of $G$ is in exactly one copy of $H_i$, for both $i=1,2$, then $G$
is said to be {\it fastened}.
If min$(m_1,m_2)=m_1=2$ and $H_1$ is a complete graph, then $G$ is said to be {\it
line-graphical}.
For example, the line graph of the $d$-cube, where
$3\le d\in\ZZ$, is a line-graphical fastened $\{K_d,K_{2,2}\}$-ultrahomogeneous
$\{K_d\}_{2^d}^2\{K_{2,2}\}_{d(d-1)2^{d-3}}^{d-1}$-graph. The first
case here, known as the cuboctahedron, is a
fastened $\{K_3,K_{2,2},C_6\}$-ultrahomogeneous
$\{K_3\}_8^2\{K_{2,2}\}_6^2\{C_6\}_4^2$-graph, where $C_6$ is 6-cycle.

In Sections 3-5, a 12-regular fastened
$\{K_4,K_{2,2,2}\}$-ultrahomogeneous
$\{K_4\}_{42}^4\{K_{2,2,2}\}_{21}^3$-graph $G$ of order 42 and
diameter 3 is presented. The role that $d$-cliques $K_d$ and squares
$K_{2,2}$ play in the line graph of the $d$-cube is performed in $G$
by tetrahedra $K_4$ and octahedra $K_{2,2,2}$, but in this case with
min$(m_1,m_2)=$min$(4,3)>2$, so $G$ is non-line-graphical.

The graph $G$ has automorphism-group order $|{\mathcal
A}(G)|=1008=4|E(G)|.$ In Section 5, the 252 edges of $G$ can be seen
as the left cosets of a subgroup $\Gamma\subset{\mathcal A}(G)$ of order
4, and its vertices as the left cosets of a subgroup of ${\mathcal
A}(G)$ of order 24.

These two equivalence classes of subgraphs of $G$, i.e. tetrahedra and
octahedra, allow in
Section 6 to define several combinatorial {\it configurations} (\cite{Cox})
related to $G$, 3 of which are {\it self-dual}, with their {\it
Levi graphs} as: {\bf(1)} a 4-regular 2-{\it arc-transitive} graph (\cite{biggs}) on 84 vertices and
1008 automorphisms, with diameter $=$ girth $=$ 6, reflecting a
natural duality property of $G$; {\bf(2)} an 8-regular
arc-transitive graph on 42 vertices and 2016 automorphisms, with
diameter $=$ 3 and girth $=$ 4; and {\bf(3)} a 6-regular {\it
semisymmetric} graph (\cite{F}) on 336 vertices and 1008
automorphisms, with diameter $=$ girth $=$ 6 and just two slightly
differing distance distributions. The {\it Menger graph} and {\it
dual Menger graph} associated to this Levi graph  have common degree
24 and diameter $=$ girth $=$ 3, with 1008 and 2016 automorphisms,
respectively.

Section 7 distinguishes the $k$-holes (or chordless $k$-cycles) of $G$ with the least $k>4$, namely $k=6$, and studies their participation in some toroidal subgraphs of $G$ that together with the octahedra of $G$ can be filled up to form a closed piecewise linear 3-manifold.

After some considerations on the Fano plane, we pass to define $G$
and study its properties.

\section{Ordered Fano pencils}

The {\it Fano plane} ${\mathcal F}$ is the $(7_3)$-configuration with points
1,2,3,4,5,6,7 and {\it Fano lines} 123, 145, 167, 246, 257, 347, 356.
The map $\Phi$ that sends the points $1,2,3,4,5,6,7$ respectively
onto the lines 123, 145, 167, 246, 257, 347, 356 has the following duality
properties: {\bf(1)} each point $p$ of ${\mathcal F}$ pertains to the lines
$\Phi(q)$, where $q\in\Phi(p)$; {\bf(2)} each Fano line
$\ell$ contains the points $\Phi(k)$, where $k$ runs over the lines
passing through $\Phi(\ell)$.

Given a point $p$ of ${\mathcal F}$, the collection of lines through $p$
is a {\it pencil} of ${\mathcal F}$. A linearly ordered presentation of
these lines is an {\it ordered pencil through} $p$. An ordered
pencil $v$ through $p$ is denoted $v=(p,q_ar_a,q_br_b,q_cr_c)$,
orderly composed, in reality, by the lines
$pq_ar_a,pq_br_b,pq_cr_c$. Note that there are $3!=6$ ordered
pencils through any point $p$ of ${\mathcal F}$.

\section{The $\{K_4,K_{2,2,2}\}$-ultrahomogeneous graph $G$}

Ordered pencils constitute
the vertex set of our claimed graph $G$, with any two vertices
$v=(p,q_ar_a,q_br_b,q_cr_c)$ and $v'=(p',q'_ar'_a,q'_br'_b,q'_cr'_c)$
adjacent whenever the following two conditions hold:
{\bf(1)} $p\ne p'$;
{\bf(2)} $|p_ir_i\cap p'_ir'_i|=1$, for $i=a,b,c$.
The 3 points of intersection resulting from item (2) form a Fano
line, which we consider as an {\it ordered Fano line} by
taking into account the subindex order $a<b<c$, and as such, set it as
the {\it strong color} of the edge $vv'$. This provides $G$ with an edge-coloring.

An alternate definition of $G$ can be
given via $\Phi^{-1}$, in which the vertices of $G$ can be seen as
the ordered Fano lines $x_ax_bx_c$, with any two such vertices
adjacent if their associated Fano lines share the entry in ${\mathcal
F}$ of exactly one of its 3 positions, either $a$ or $b$ or $c$.
We keep throughout, however, the ordered-pencil presentation of $G$,
but the first self-dual configuration of Subsection 6.1 and
accompanying example show that the suggested dual presentation of
$G$ is valid as well.

Notice that the vertices of $G$ with initial entry $p=1$ appear in lexicographic order as:
$$(1,23,45,67),(1,23,67,45),(1,45,23,67),(1,45,67,23),(1,67,23,45),(1,67,45,23),$$
which may be simplified in notation by using super-indices $a$ through $f$ to
denote the shown order, that is:
$1^a,1^b,1^c,1^d,1^e,1^f,$ respectively.
A similar lexicographic presentation may be given to the vertices
of $G$ having $p=2,\ldots,7$. This treatment covers the 42 vertices
of $G$. As an example of the adjacency of $G$,
the neighbors of $1^a=(1,23,45,67)$ in $G$ are:
$$\begin{array}{cccc}
(2,13,46,57), & (2,13,57,46), & (3,12,47,56), & (3,12,56,47), \\
(4,26,15,37), & (4,37,15,26), & (5,27,14,36), & (5,36,14,27), \\
(6,24,35,17), & (6,35,24,17), & (7,25,34,16), & (7,34,25,16),
\end{array}$$
or in the continuation of the simplified notation above:
$2^a,2^b,3^a,3^b,4^c,4^e,5^c,5^e,6^d,6^f,7^d,7^f$.
The strong colors of the resulting edges are:
$167$, $154$, $176$, $154$, $356$, $246$, $347$, $451$, $321$, $231$, $321$, respectively.

Given vertices $v=(p,q_ar_a,q_br_b,q_cr_c)$ and
$w=(p',q'_ar'_a,q'_br'_b,q'_cr'_c)$ adjacent in $G$, there exists a
well-defined $j\in\{a,b,c\}$ such that {\bf(1)} $p\in q'_jr'_j$;
{\bf(2)} $p'\in q_jr_j$; {\bf(3)} the lines $pq_jr_j$ and
$p'q'_jr'_j$ intersect at either $q_j$ or $r_j$, which coincides
with either $q'_j$ or $r'_j$. Say that these lines $pq_jr_j$ and
$p'q'_jr'_j$ intersect at $q_j$. Then $q_j$ (including the subindex
$j$) is taken as the {\it weak color} for the edge $vw$. This
provides $G$ with another edge-coloring, with symbols $q_j$, where $q\in {\mathcal
F}$ and $j\in\{a,b,c\}$. For example, the weak colors $q_j$
corresponding to the 12 edges incident to $1^a$, as cited above,
are: $3_a,$ $3_a,$ $2_a,$ $2_a,$ $5_b,$ $5_b,$ $4_b,$ $4_b,$ $7_c,$
$7_c,$ $6_c,$ $6_c$, respectively.

\subsection{The automorphism group ${\mathcal A}(G)$ of $G$}

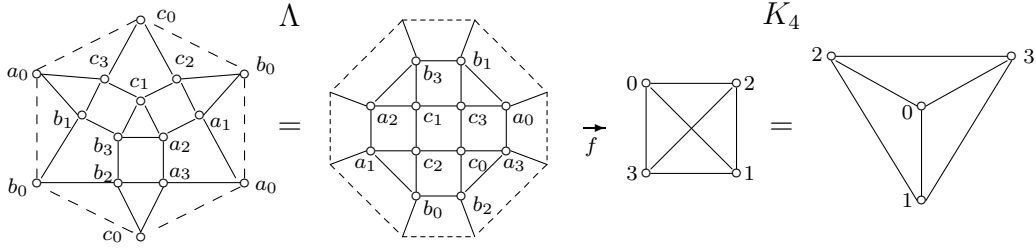
\begin{figure}[htp]
\unitlength=0.60mm
\special{em:linewidth 0.4pt}
\linethickness{0.4pt}
\begin{picture}(249.00,50.00)
\put(44.00,21.00){\makebox(0,0)[cc]{$_{b_3}$}}
\put(25.00,37.00){\makebox(0,0)[cc]{$_{a_0}$}}
\put(47.00,23.00){\circle{2.00}}
\put(57.00,23.00){\circle{2.00}}
\put(52.00,31.00){\circle{2.00}}
\emline{48.00}{23.00}{1}{56.00}{23.00}{2}
\put(65.00,28.00){\circle{2.00}}
\put(60.00,36.00){\circle{2.00}}
\emline{56.00}{24.00}{3}{53.00}{30.00}{4}
\emline{53.00}{32.00}{5}{59.00}{35.00}{6}
\emline{61.00}{35.00}{7}{64.00}{29.00}{8}
\emline{64.00}{27.00}{9}{58.00}{24.00}{10}
\put(39.00,28.00){\circle{2.00}}
\put(44.00,36.00){\circle{2.00}}
\emline{46.00}{24.00}{11}{40.00}{27.00}{12}
\emline{40.00}{29.00}{13}{43.00}{35.00}{14}
\emline{45.00}{35.00}{15}{51.00}{32.00}{16}
\emline{48.00}{24.00}{17}{51.00}{30.00}{18}
\put(47.00,13.00){\circle{2.00}}
\put(57.00,13.00){\circle{2.00}}
\emline{47.00}{22.00}{19}{47.00}{14.00}{20}
\emline{48.00}{13.00}{21}{56.00}{13.00}{22}
\emline{57.00}{14.00}{23}{57.00}{22.00}{24}
\put(29.00,13.00){\circle{2.00}}
\put(75.00,13.00){\circle{2.00}}
\put(80.00,37.00){\makebox(0,0)[cb]{$_{b_0}$}}
\put(29.00,37.00){\circle{2.00}}
\put(75.00,37.00){\circle{2.00}}
\put(52.00,1.00){\circle{2.00}}
\put(52.00,49.00){\circle{2.00}}
\emline{33.00}{40.00}{25}{35.00}{41.00}{26}
\emline{37.00}{42.00}{27}{39.00}{43.00}{28}
\emline{41.00}{44.00}{29}{43.00}{45.00}{30}
\emline{45.00}{46.00}{31}{47.00}{47.00}{32}
\emline{49.00}{48.00}{33}{51.00}{49.00}{34}
\emline{56.00}{3.00}{35}{58.00}{4.00}{36}
\emline{60.00}{5.00}{37}{62.00}{6.00}{38}
\emline{64.00}{7.00}{39}{66.00}{8.00}{40}
\emline{68.00}{9.00}{41}{70.00}{10.00}{42}
\emline{72.00}{11.00}{43}{74.00}{12.00}{44}
\emline{50.00}{2.00}{45}{50.00}{2.00}{46}
\emline{50.00}{2.00}{47}{48.00}{3.00}{48}
\emline{46.00}{4.00}{49}{44.00}{5.00}{50}
\emline{42.00}{6.00}{51}{40.00}{7.00}{52}
\emline{38.00}{8.00}{53}{36.00}{9.00}{54}
\emline{34.00}{10.00}{55}{32.00}{11.00}{56}
\emline{73.00}{39.00}{57}{73.00}{39.00}{58}
\emline{73.00}{39.00}{59}{71.00}{40.00}{60}
\emline{69.00}{41.00}{61}{67.00}{42.00}{62}
\emline{65.00}{43.00}{63}{63.00}{44.00}{64}
\emline{61.00}{45.00}{65}{59.00}{46.00}{66}
\emline{57.00}{47.00}{67}{55.00}{48.00}{68}
\emline{75.00}{36.00}{69}{75.00}{34.00}{70}
\emline{75.00}{32.00}{71}{75.00}{30.00}{72}
\emline{75.00}{28.00}{73}{75.00}{26.00}{74}
\emline{75.00}{24.00}{75}{75.00}{22.00}{76}
\emline{75.00}{20.00}{77}{75.00}{18.00}{78}
\emline{75.00}{16.00}{79}{75.00}{14.00}{80}
\emline{29.00}{36.00}{81}{29.00}{34.00}{82}
\emline{29.00}{32.00}{83}{29.00}{30.00}{84}
\emline{29.00}{28.00}{85}{29.00}{26.00}{86}
\emline{29.00}{24.00}{87}{29.00}{22.00}{88}
\emline{29.00}{20.00}{89}{29.00}{18.00}{90}
\emline{29.00}{16.00}{91}{29.00}{14.00}{92}
\emline{74.00}{13.00}{93}{58.00}{13.00}{94}
\emline{46.00}{13.00}{95}{30.00}{13.00}{96}
\emline{38.00}{28.00}{97}{30.00}{37.00}{98}
\emline{61.00}{36.00}{99}{74.00}{37.00}{100}
\emline{74.00}{37.00}{101}{66.00}{28.00}{102}
\emline{43.00}{36.00}{103}{30.00}{37.00}{104}
\emline{47.00}{12.00}{105}{52.00}{2.00}{106}
\emline{52.00}{2.00}{107}{57.00}{12.00}{108}
\put(80.00,12.00){\makebox(0,0)[cc]{$_{a_0}$}}
\put(58.00,50.00){\makebox(0,0)[cc]{$_{c_0}$}}
\put(70.00,26.00){\makebox(0,0)[cc]{$_{a_1}$}}
\put(62.00,40.00){\makebox(0,0)[cc]{$_{c_2}$}}
\put(52.00,35.00){\makebox(0,0)[cc]{$_{c_1}$}}
\put(61.00,21.00){\makebox(0,0)[cc]{$_{a_2}$}}
\put(43.00,40.00){\makebox(0,0)[cc]{$_{c_3}$}}
\put(61.00,15.00){\makebox(0,0)[cc]{$_{a_3}$}}
\put(35.00,27.00){\makebox(0,0)[cc]{$_{b_1}$}}
\put(44.00,15.00){\makebox(0,0)[cc]{$_{b_2}$}}
\put(25.00,12.00){\makebox(0,0)[cc]{$_{b_0}$}}
\put(46.00,1.00){\makebox(0,0)[cc]{$_{c_0}$}}
\emline{29.00}{38.00}{109}{31.00}{39.00}{110}
\emline{53.00}{48.00}{111}{59.00}{37.00}{112}
\emline{66.00}{27.00}{113}{73.00}{14.00}{114}
\emline{51.00}{48.00}{115}{45.00}{37.00}{116}
\emline{38.00}{27.00}{117}{30.00}{14.00}{118}
\emline{110.00}{49.00}{119}{111.00}{49.00}{120}
\emline{112.00}{49.00}{121}{113.00}{49.00}{122}
\emline{114.00}{49.00}{123}{115.00}{49.00}{124}
\emline{116.00}{49.00}{125}{117.00}{49.00}{126}
\emline{118.00}{49.00}{127}{119.00}{49.00}{128}
\emline{120.00}{49.00}{129}{121.00}{49.00}{130}
\emline{122.00}{49.00}{131}{123.00}{49.00}{132}
\emline{124.00}{49.00}{133}{125.00}{49.00}{134}
\emline{94.00}{33.00}{135}{94.00}{32.00}{136}
\emline{94.00}{31.00}{137}{94.00}{30.00}{138}
\emline{94.00}{29.00}{139}{94.00}{28.00}{140}
\emline{94.00}{27.00}{141}{94.00}{26.00}{142}
\emline{94.00}{25.00}{143}{94.00}{24.00}{144}
\emline{94.00}{23.00}{145}{94.00}{22.00}{146}
\emline{94.00}{21.00}{147}{94.00}{20.00}{148}
\emline{94.00}{19.00}{149}{94.00}{18.00}{150}
\emline{94.00}{17.00}{151}{95.00}{16.00}{152}
\emline{96.00}{15.00}{153}{97.00}{14.00}{154}
\emline{98.00}{13.00}{155}{99.00}{12.00}{156}
\emline{100.00}{11.00}{157}{101.00}{10.00}{158}
\emline{102.00}{9.00}{159}{103.00}{8.00}{160}
\emline{104.00}{7.00}{161}{105.00}{6.00}{162}
\emline{106.00}{5.00}{163}{107.00}{4.00}{164}
\emline{108.00}{3.00}{165}{109.00}{2.00}{166}
\emline{95.00}{34.00}{167}{96.00}{35.00}{168}
\emline{97.00}{36.00}{169}{98.00}{37.00}{170}
\emline{99.00}{38.00}{171}{100.00}{39.00}{172}
\emline{101.00}{40.00}{173}{102.00}{41.00}{174}
\emline{103.00}{42.00}{175}{104.00}{43.00}{176}
\emline{105.00}{44.00}{177}{106.00}{45.00}{178}
\emline{107.00}{46.00}{179}{108.00}{47.00}{180}
\emline{109.00}{48.00}{181}{110.00}{49.00}{182}
\put(164.00,35.00){\circle{2.00}}
\put(184.00,35.00){\circle{2.00}}
\put(164.00,15.00){\circle{2.00}}
\put(184.00,15.00){\circle{2.00}}
\emline{165.00}{15.00}{183}{183.00}{15.00}{184}
\emline{165.00}{35.00}{185}{183.00}{35.00}{186}
\emline{164.00}{34.00}{187}{164.00}{16.00}{188}
\emline{184.00}{16.00}{189}{184.00}{34.00}{190}
\put(117.00,27.00){\makebox(0,0)[cc]{$_{c_1}$}}
\put(117.00,17.00){\makebox(0,0)[cc]{$_{c_2}$}}
\put(127.00,27.00){\makebox(0,0)[cc]{$_{c_3}$}}
\put(127.00,17.00){\makebox(0,0)[cc]{$_{c_0}$}}
\put(137.00,27.00){\makebox(0,0)[cc]{$_{a_0}$}}
\put(135.00,17.00){\makebox(0,0)[cc]{$_{a_3}$}}
\put(107.00,27.00){\makebox(0,0)[cc]{$_{a_2}$}}
\put(102.00,17.00){\makebox(0,0)[cc]{$_{a_1}$}}
\put(117.00,37.00){\makebox(0,0)[cc]{$_{b_3}$}}
\put(128.00,40.00){\makebox(0,0)[cc]{$_{b_1}$}}
\put(117.00,7.00){\makebox(0,0)[cc]{$_{b_0}$}}
\put(128.00,8.00){\makebox(0,0)[cc]{$_{b_2}$}}
\put(161.00,35.00){\makebox(0,0)[cc]{$_0$}}
\put(161.00,15.00){\makebox(0,0)[cc]{$_3$}}
\put(187.00,35.00){\makebox(0,0)[cc]{$_2$}}
\put(187.00,15.00){\makebox(0,0)[cc]{$_1$}}
\emline{165.00}{34.00}{191}{183.00}{16.00}{192}
\emline{183.00}{34.00}{193}{165.00}{16.00}{194}
\put(150.00,25.00){\vector(1,0){5.00}}
\put(152.00,21.00){\makebox(0,0)[cc]{$_f$}}
\put(85.00,50.00){\makebox(0,0)[cc]{$\Lambda$}}
\put(194.00,50.00){\makebox(0,0)[cc]{$K_4$}}
\put(85.00,25.00){\makebox(0,0)[cc]{$=$}}
\put(194.00,25.00){\makebox(0,0)[cc]{$=$}}
\put(222.00,9.00){\makebox(0,0)[cc]{$_1$}}
\put(225.00,9.00){\circle{2.00}}
\put(225.00,30.00){\circle{2.00}}
\emline{225.00}{29.00}{195}{225.00}{10.00}{196}
\emline{244.00}{40.00}{197}{226.00}{30.00}{198}
\emline{206.00}{40.00}{199}{224.00}{30.00}{200}
\put(222.00,27.00){\makebox(0,0)[cb]{$_0$}}
\emline{245.00}{40.00}{201}{226.00}{9.00}{202}
\emline{224.00}{9.00}{203}{205.00}{40.00}{204}
\put(202.00,40.00){\makebox(0,0)[cb]{$_2$}}
\put(249.00,41.00){\makebox(0,0)[cc]{$_3$}}
\put(205.00,41.00){\circle{2.00}}
\put(245.00,41.00){\circle{2.00}}
\emline{206.00}{41.00}{205}{244.00}{41.00}{206}
\put(113.00,30.00){\circle{2.00}}
\put(123.00,30.00){\circle{2.00}}
\emline{114.00}{30.00}{207}{122.00}{30.00}{208}
\put(123.00,20.00){\circle{2.00}}
\emline{113.00}{29.00}{209}{113.00}{21.00}{210}
\emline{114.00}{20.00}{211}{122.00}{20.00}{212}
\emline{123.00}{21.00}{213}{123.00}{29.00}{214}
\put(113.00,40.00){\circle{2.00}}
\emline{113.00}{39.00}{215}{113.00}{31.00}{216}
\put(113.00,20.00){\circle{2.00}}
\emline{113.00}{19.00}{217}{113.00}{11.00}{218}
\put(113.00,10.00){\circle{2.00}}
\put(123.00,40.00){\circle{2.00}}
\emline{123.00}{39.00}{219}{123.00}{31.00}{220}
\emline{123.00}{19.00}{221}{123.00}{11.00}{222}
\put(123.00,10.00){\circle{2.00}}
\put(103.00,30.00){\circle{2.00}}
\emline{104.00}{30.00}{223}{112.00}{30.00}{224}
\put(133.00,30.00){\circle{2.00}}
\emline{124.00}{30.00}{225}{132.00}{30.00}{226}
\put(103.00,20.00){\circle{2.00}}
\emline{104.00}{20.00}{227}{112.00}{20.00}{228}
\put(133.00,20.00){\circle{2.00}}
\emline{124.00}{20.00}{229}{132.00}{20.00}{230}
\emline{104.00}{31.00}{231}{112.00}{39.00}{232}
\emline{124.00}{39.00}{233}{132.00}{31.00}{234}
\emline{132.00}{19.00}{235}{124.00}{11.00}{236}
\emline{112.00}{11.00}{237}{104.00}{19.00}{238}
\emline{103.00}{21.00}{239}{103.00}{29.00}{240}
\emline{114.00}{40.00}{241}{122.00}{40.00}{242}
\emline{133.00}{29.00}{243}{133.00}{21.00}{244}
\emline{122.00}{10.00}{245}{114.00}{10.00}{246}
\emline{123.00}{41.00}{247}{126.00}{49.00}{248}
\emline{113.00}{41.00}{249}{110.00}{49.00}{250}
\emline{102.00}{30.00}{251}{94.00}{33.00}{252}
\emline{102.00}{20.00}{253}{94.00}{17.00}{254}
\emline{110.00}{1.00}{255}{113.00}{9.00}{256}
\emline{126.00}{1.00}{257}{123.00}{9.00}{258}
\emline{110.00}{1.00}{259}{111.00}{1.00}{260}
\emline{112.00}{1.00}{261}{113.00}{1.00}{262}
\emline{114.00}{1.00}{263}{115.00}{1.00}{264}
\emline{116.00}{1.00}{265}{117.00}{1.00}{266}
\emline{118.00}{1.00}{267}{119.00}{1.00}{268}
\emline{120.00}{1.00}{269}{121.00}{1.00}{270}
\emline{122.00}{1.00}{271}{123.00}{1.00}{272}
\emline{124.00}{1.00}{273}{125.00}{1.00}{274}
\emline{142.00}{17.00}{275}{134.00}{20.00}{276}
\emline{142.00}{33.00}{277}{134.00}{30.00}{278}
\emline{126.00}{49.00}{279}{127.00}{48.00}{280}
\emline{128.00}{47.00}{281}{129.00}{46.00}{282}
\emline{130.00}{45.00}{283}{131.00}{44.00}{284}
\emline{132.00}{43.00}{285}{133.00}{42.00}{286}
\emline{134.00}{41.00}{287}{135.00}{40.00}{288}
\emline{136.00}{39.00}{289}{137.00}{38.00}{290}
\emline{138.00}{37.00}{291}{139.00}{36.00}{292}
\emline{140.00}{35.00}{293}{141.00}{34.00}{294}
\emline{142.00}{33.00}{295}{142.00}{32.00}{296}
\emline{142.00}{31.00}{297}{142.00}{30.00}{298}
\emline{142.00}{29.00}{299}{142.00}{28.00}{300}
\emline{142.00}{27.00}{301}{142.00}{26.00}{302}
\emline{142.00}{25.00}{303}{142.00}{24.00}{304}
\emline{142.00}{23.00}{305}{142.00}{22.00}{306}
\emline{142.00}{21.00}{307}{142.00}{20.00}{308}
\emline{142.00}{19.00}{309}{142.00}{18.00}{310}
\emline{127.00}{2.00}{311}{128.00}{3.00}{312}
\emline{129.00}{4.00}{313}{130.00}{5.00}{314}
\emline{131.00}{6.00}{315}{132.00}{7.00}{316}
\emline{133.00}{8.00}{317}{134.00}{9.00}{318}
\emline{135.00}{10.00}{319}{136.00}{11.00}{320}
\emline{137.00}{12.00}{321}{138.00}{13.00}{322}
\emline{139.00}{14.00}{323}{140.00}{15.00}{324}
\emline{141.00}{16.00}{325}{142.00}{17.00}{326}
\end{picture}
\caption{The homomorphism $f:\Lambda\rightarrow K_4$
with $f(j_i)=i$ for $i\in\{0,1,2,3\}$, $j\in\{a,b,c\}$}
\end{figure}

The 12 neighbors of $1^a$ displayed above induce a subgraph $N_G(1^a)$ of $G$, called the
open neighborhood of $1^a$ in $G$, which is isomorphic to the graph $\Lambda$ of the
hemi-rhom\-bi\-cu\-boc\-ta\-he\-dron (obtained from the rhombicuboctahedron by identification
of antipodal vertices and edges). This is a 4-regular vertex-transitive graph on 12 vertices
embedded in the projective plane with 13 faces realized by
4 disjoint triangles and 9 additional
4-holes. The 4-holes are of two types: {\bf(1)}
6 have two opposite sides adjacent each to a triangle; {\bf(2)}
the other 3 have only their vertices in common with the 4 triangles. We also have the
graph homomorphism $f:\Lambda\rightarrow K_4$ of
Figure 1, where
$f(j_i)=i$ for $i\in\{0,1,2,3\}$, $j\in\{a,b,c\}$ and
$\Lambda$ is depicted in two different ways inside (dotted) fundamental polygons of the real projective plane. Moreover, we may identify
$\Lambda$ with $N_G(1^a)$ via a graph isomorphism $g:\Lambda\rightarrow N_G(1^a)$ given by:

$$\begin{array}{cccc}
g(a_0)=5^c, & g(a_1)=4^c, & g(a_2)=5^e, & g(a_3)=4^e, \\
g(b_0)=6^d, & g(b_1)=7^d, & g(b_2)=7^f, & g(b_3)=6^f, \\
g(c_0)=2^b, & g(c_1)=2^a, & g(c_2)=3^b, & g(c_3)=3^a.
\end{array}$$
Moreover, the graph homomorphism $f$ induces, at the level of automorphism groups of graphs, a group isomorphism $f^*:{\mathcal A}(\Lambda)\rightarrow{\mathcal A}(K_4)=S_4$. In fact, $f^*$ is given by sending the following generators of ${\mathcal A}(\Lambda)$ into corresponding generators of $S_4$ (that can be better visualized from the leftmost $\Lambda$ to the rightmost $K_4$ depicted in Figure 1):
$$\begin{array}{lll}
(a_1b_2c_3)(b_1c_2a_3)(c_1a_2b_3)(a_0b_0c_0) & \rightarrow & (123), \\
(a_0a_1)(a_2a_3)(b_0b_1)(b_2b_3)(c_0c_1)(c_2c_3) & \rightarrow & (01)(23), \\
(a_0a_1a_2a_3)(c_0b_1c_2b_3)(b_0c_1b_2c_3) & \rightarrow & (0123).
\end{array}$$
Thus, ${\mathcal A}(N_G(1^a))={\mathcal A}(\Lambda)=S_4$ has 24 elements, which is consistent
with the size of a vertex stabilizer of $G$. Furthermore, since $G$ has 42 vertices that behave exactly in the same geometric way as ordered pencils in $\mathcal F$, we conclude that $|{\mathcal A}(G)|=42\times 24=1008$.

\section{Copies of $K_{2,2,2}$ and $K_4$ in $G$}

Notice that $f$ maps bijectively the 3 4-cycles and 4
triangles of $K_4$ respectively onto the 3 4-holes of $\Lambda$ of type (2)
above and the 4 triangles of $\Lambda$. Notice also that these 7 holes form a cycle-decomposition of $\Lambda$. Inside the closed neighborhood $N_G[w]$ of each vertex $w$ of $G$ (induced in $G$ by $w$ and the open neighborhood $N_G(w)$),
we obtain 3 copies of $K_{2,2,2}$ and 4 copies of $K_4$, which are induced by $w$ together respectively with the mentioned 3 4-holes and 4 triangles. Observe that these 7 induced subgraphs of $G$ have intersection formed solely  by $w$. The rest of this section is dedicated to the study of these polyhedral subgraphs.

First, notice that the inverse image $f^{-1}$ of each edge of $K_4$ is one of the 6 4-holes of
$\Lambda$ of type (1) above.
This yields another cycle-decomposition of $\Lambda$, which in turn makes explicit the remaining 4-holes of $G$, apart from the 4-holes contained in the copies of $K_{2,2,2}$ of $G$. However, these new 4-holes are not contained in any copy of $K_{2,2,2}$ in $G$.

\subsection{Copies of $K_{2,2,2}$ in $G$}

Each vertex of $G$ belongs to 3 induced copies of $K_{2,2,2}$ in
$G$. For example, the sets of weak colors $q_j$ of the edges of such copies for the vertex $1^a$, that contain the 4-holes $g(c_0,c_2,c_1,c_3)$, $g(a_0,a_1,a_2,a_3)$ and $g(b_0,b_1,b_3,b_2)$ arising in Subsection 3.1, are respectively: $\{1_a,2_a,3_a\},\{1_b,4_b,5_b\},\{1_c,6_c,7_c\}$.

Each $q_j$ colors
the edges of a specific 4-hole in its corresponding copy of $K_{2,2,2}$. The 3 weak colors appearing in each copy of $K_{2,2,2}$ correspond bijectively with its 3
4-holes, the edges of each 4-hole bearing a common weak color of its own.

A similar situation holds for any other vertex of $G$. There is a
copy of $K_{2,2,2}$ in $G$ whose set of weak colors of edges is
$\{x_j,y_j,z_j\}$, for each line $xyz$ of ${\mathcal F}$ and index
$j\in\{a,b,c\}$. We denote this copy of $K_{2,2,2}$ by $[xyz]_j$. As
a result, there is a total of $21=7\times 3$ copies of $K_{2,2,2}$
in $G$. In fact, triangles with weak colors $q_j$ sharing a common
$j$ (but $q$ varying) are
only present in the said copies of $K_{2,2,2}$ in $G$. Each
4-hole in a copy of $K_{2,2,2}$ in $G$ have: {\bf(1)} edges sharing a common weak color $q_j$ and {\bf(2)} opposite vertices representing ordered pencils through a common
point of ${\mathcal F}$, which yields a total of two such points per
4-hole.

For example, the strong colors of the triangles $[xyz]_j$ composing the copies of $K_{2,2,2}$
incident to $1_a$ conform triples of strong colors having:

\noindent for $[123]_a$, $a$-entries covering line 123, and another
fixed entry equal to each one of 4,5,6,7:
$$\begin{array}{cccc}
(145,246,347), & (154,257,356), & (167,257,347), & (176,275,374), \\
(154,264,374), & (145,275,365), & (176,275,374), & (167,246,365);
\end{array}$$
for $[145]_b$, $b$-entries covering line 145, and another fixed entry
equal to each one of 2,3,6,7:
$$\begin{array}{cccc}
(213,246,257), & (312,347,356), & (617,642,653), & (716,743,752), \\
(312,642,752), & (213,743,653), & (716,246,356), & (617,347,257);
\end{array}$$
for $[167]_c$, $c$-entries covering line 167, and another fixed entry
equal to each one of 2,3,4,5:
$$\begin{array}{cccc}
(231,246,257), & (321,356,347), & (451,426,437), & (541,536,527), \\
(321,426,527), & (231,536,437), & (541,246,347), & (451,356,257).
\end{array}$$
In fact, these triangles are respectively:
$$\begin{array}{l}
(2^a,3^a,1^a),(2^b,3^b,1^a),(2^a,3^b,1^a),(2^b,3^a,1^a),(2^b,3^b,1^b),(2^a,3^a,1^b),(2^b,3^a,1^b),(2^a,3^b,1^b); \\
(4^c,5^c,1^a),(4^e,5^e,1^a),(4^c,5^e,1^f),(4^e,5^c,1^f),(4^e,5^e,1^f),(4^c,5^c,1^f),(4^e,5^c,1^a),(4^c,5^e,1^a); \\
(6^d,7^d,\!1^a),(6^f,7^f,\!1^a),(6^d,7^f,1^c),(6^f,7^d,\!1^c),(6^f,7^f,\!1^c),(6^d,7^d,1^c),(6^f,7^d,1^a),(6^d,7^f,1^a).
\end{array}$$
The sets of strong colors for the respective composing 4-holes are:
$$\begin{array}{ccc}
(145,167,154,176), & (246,257,264,275), & (347,356,374,365); \\
(213,617,312,716), & (246,347,624,743), & (257,356,752,653); \\
(231,451,321,541), & (246,356,426,536), & (257,347,527,437).
\end{array}$$
In fact, these 4-holes are respectively:
$$\begin{array}{l}
(2^a,3^a,2^b,3^b),(3^a,1^a,3^b,1^b),(2^a,1^a,2^b,1^b); \\
(4^c,5^e,4^e,5^c),(5^c,1^a,5^e,1^f),(4^c,1^a,4^e,1^f); \\
(7^d,6^d,7^f,\!6^f),(7^d,1^a,7^f,\!1^c),(6^d,1^a,6^f,1^c).
\end{array}$$

\subsection{Copies of $K_4$ in $G$}

There is one copy of $K_4$ in $G$ for each ordered Fano line $xyz$.
Such a copy, denoted $\langle xyz\rangle $, is formed by 3 pairs of equally weakly-colored
opposite edges, with weak colors $x_a$, $y_b$ and $z_c$.
For each $p\in {\mathcal F}\setminus\{x,y,z\}$, there is exactly one
vertex $(p,q_ar_a,q_br_b,q_cr_c)$ of $\langle xyz\rangle $, with
$x\in q_ar_a$, $y\in q_br_b$, $z\in q_cr_c$. The strong colors of the edges of
$\langle xyz\rangle $ are precisely $xyz$.
For example, the triangles $g(c_1,a_2,b_2)$, $g(c_3,a_0,b_1)$, $g(c_2,a_1,b_0)$ and $g(c_0,a_3,b_2)$ from Subsection 3.1 are contained respectively in $\langle 347\rangle$,
$\langle 246\rangle$, $\langle 257\rangle$ and $\langle 356\rangle$.
Since there are 42 such copies of $K_4$ in $G$, we arrive at the following result.

\begin{teo} The graph $G$ is a 12-regular
$\{K_4\}_{42}^4\{K_{2,2,2}\}_{21}^3$-graph of order 42 and
diameter 3. Each vertex of $G$ is incident to exactly 3
copies of $K_{2,2,2}$ and 4 copies of $K_4$.
\end{teo}

\begin{demo} Let $G'$ be the graph defined by the same rules that define $G$ on the unordered Fano lines. Then it is not hard to prove that $G'$ is isomorphic to the graph $2K_7$, the complete graph on 7 vertices with each edge doubled. The graph $G$ is then a 6-fold covering graph over $G'$. Also, the lexicographically smallest path
realizing the diameter of $G$ is the 3-path $(1^a$, $2^a$, $4^a$, $1^d).$
The statement follows.
\end{demo}

\subsection{Disposition of copies of $K_{2,2,2}$ and $K_4$ in $G$}

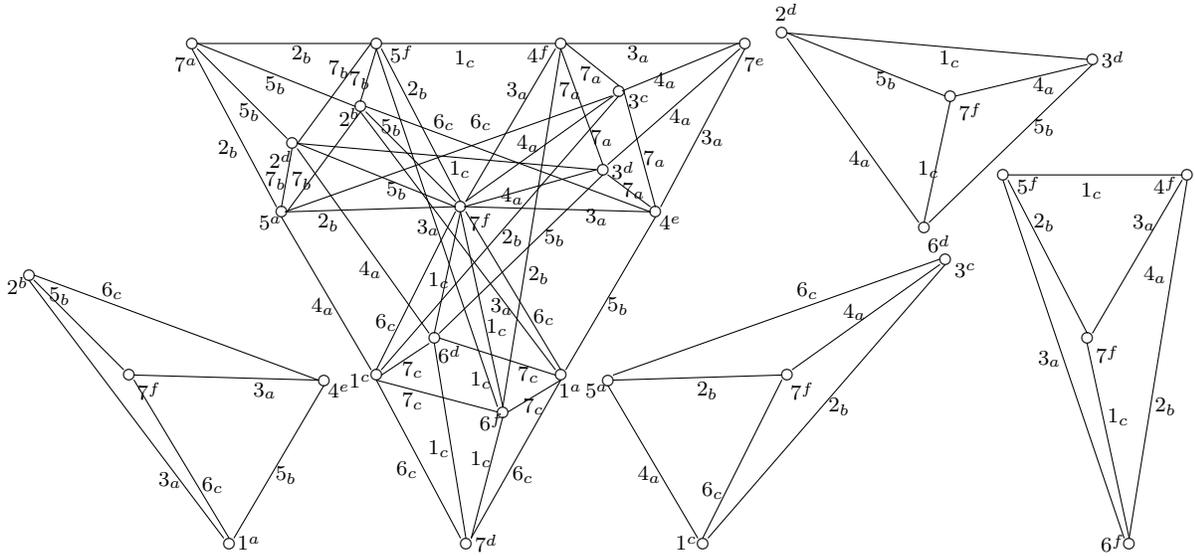
\begin{figure}[htp]
\unitlength=0.70mm
\special{em:linewidth 0.4pt} \linethickness{0.4pt}
\begin{picture}(225.33,102.00)
\put(35.33,96.00){\circle{2.00}} \put(70.33,96.00){\circle{2.00}}
\put(105.33,96.00){\circle{2.00}} \put(140.33,96.00){\circle{2.00}}
\emline{36.33}{96.00}{1}{69.33}{96.00}{2}
\emline{71.33}{96.00}{3}{104.33}{96.00}{4}
\emline{106.33}{96.00}{5}{139.33}{96.00}{6}
\put(87.33,1.00){\circle{2.00}} \put(70.33,33.00){\circle{2.00}}
\put(123.33,64.00){\circle{2.00}}
\emline{70.33}{95.00}{7}{70.33}{95.00}{8}
\put(67.33,84.00){\circle{2.00}} \put(116.33,87.00){\circle{2.00}}
\put(113.33,72.00){\circle{2.00}} \put(94.33,26.00){\circle{2.00}}
\put(81.33,40.00){\circle{2.00}}
\emline{35.33}{95.00}{9}{51.33}{65.00}{10}
\emline{52.33}{63.00}{11}{69.33}{34.00}{12}
\emline{70.33}{32.00}{13}{86.33}{2.00}{14}
\emline{88.33}{2.00}{15}{105.33}{32.00}{16}
\emline{106.33}{34.00}{17}{123.33}{63.00}{18}
\emline{124.33}{65.00}{19}{140.33}{95.00}{20}
\emline{36.33}{96.00}{21}{66.33}{84.00}{22}
\emline{68.33}{84.00}{23}{122.33}{64.00}{24}
\emline{70.33}{95.00}{25}{67.33}{85.00}{26}
\emline{67.33}{83.00}{27}{53.33}{64.00}{28}
\put(52.33,64.00){\circle{2.00}}
\emline{52.33}{65.00}{29}{54.33}{76.00}{30}
\emline{55.33}{77.00}{31}{69.33}{96.00}{32}
\put(54.33,77.00){\circle{2.00}}
\emline{106.33}{96.00}{33}{116.33}{88.00}{34}
\emline{117.33}{87.00}{35}{123.33}{65.00}{36}
\emline{122.33}{65.00}{37}{114.33}{71.00}{38}
\emline{113.33}{73.00}{39}{105.33}{95.00}{40}
\emline{82.33}{40.00}{41}{104.33}{33.00}{42}
\emline{105.33}{32.00}{43}{95.33}{26.00}{44}
\emline{93.33}{26.00}{45}{70.33}{32.00}{46}
\put(105.33,33.00){\circle{2.00}}
\emline{88.33}{2.00}{47}{94.33}{25.00}{48}
\emline{36.33}{95.00}{49}{53.33}{78.00}{50}
\emline{139.33}{95.00}{51}{114.33}{73.00}{52}
\emline{87.33}{2.00}{53}{81.33}{39.00}{54}
\emline{68.33}{83.00}{55}{85.33}{66.00}{56}
\put(86.33,65.00){\circle{2.00}}
\emline{55.33}{77.00}{57}{85.33}{65.00}{58}
\emline{53.33}{64.00}{59}{85.33}{65.00}{60}
\emline{71.33}{95.00}{61}{86.33}{66.00}{62}
\emline{87.33}{66.00}{63}{104.33}{95.00}{64}
\emline{87.33}{66.00}{65}{115.33}{86.00}{66}
\emline{87.33}{65.00}{67}{122.33}{64.00}{68}
\emline{86.33}{64.00}{69}{81.33}{41.00}{70}
\emline{85.33}{64.00}{71}{70.33}{34.00}{72}
\emline{87.33}{64.00}{73}{105.33}{34.00}{74}
\emline{86.33}{64.00}{75}{94.33}{27.00}{76}
\emline{67.33}{83.00}{77}{104.33}{34.00}{78}
\emline{71.33}{33.00}{79}{80.33}{39.00}{80}
\emline{71.33}{34.00}{81}{116.33}{86.00}{82}
\emline{70.33}{95.00}{83}{93.33}{27.00}{84}
\emline{94.33}{27.00}{85}{105.33}{95.00}{86}
\emline{55.33}{77.00}{87}{112.33}{72.00}{88}
\emline{113.33}{71.00}{89}{82.33}{41.00}{90}
\emline{80.33}{41.00}{91}{55.33}{76.00}{92}
\put(126.33,62.00){\makebox(0,0)[cc]{$_{4^e}$}}
\put(120.33,85.00){\makebox(0,0)[cc]{$_{3^c}$}}
\put(67.33,32.00){\makebox(0,0)[cc]{$_{1^c}$}}
\put(91.33,1.00){\makebox(0,0)[cc]{$_{7^d}$}}
\put(107.33,30.00){\makebox(0,0)[cc]{$_{1^a}$}}
\put(75.33,94.00){\makebox(0,0)[cc]{$_{5^f}$}}
\put(52.33,74.00){\makebox(0,0)[cc]{$_{2^d}$}}
\put(50.33,62.00){\makebox(0,0)[cc]{$_{5^a}$}}
\put(65.33,82.00){\makebox(0,0)[cc]{$_{2^b}$}}
\put(92.33,24.00){\makebox(0,0)[cc]{$_{6^f}$}}
\put(84.33,37.00){\makebox(0,0)[cc]{$_{6^d}$}}
\put(117.33,72.00){\makebox(0,0)[cc]{$_{3^d}$}}
\put(101.33,94.00){\makebox(0,0)[cc]{$_{4^f}$}}
\put(34.33,92.00){\makebox(0,0)[cc]{$_{7^a}$}}
\put(142.33,92.00){\makebox(0,0)[cc]{$_{7^e}$}}
\put(90.33,62.00){\makebox(0,0)[cc]{$_{7^f}$}}
\put(56.33,94.00){\makebox(0,0)[cc]{$_{2_b}$}}
\put(46.33,83.00){\makebox(0,0)[cc]{$_{5_b}$}}
\put(51.33,88.00){\makebox(0,0)[cc]{$_{5_b}$}}
\put(42.33,76.00){\makebox(0,0)[cc]{$_{2_b}$}}
\put(61.33,62.00){\makebox(0,0)[cc]{$_{2_b}$}}
\put(134.33,78.00){\makebox(0,0)[cc]{$_{3_a}$}}
\put(125.33,89.00){\makebox(0,0)[cc]{$_{4_a}$}}
\put(128.33,82.00){\makebox(0,0)[cc]{$_{4_a}$}}
\put(120.33,94.00){\makebox(0,0)[cc]{$_{3_a}$}}
\put(76.33,15.00){\makebox(0,0)[cc]{$_{6_c}$}}
\put(98.33,14.00){\makebox(0,0)[cc]{$_{6_c}$}}
\put(82.33,19.00){\makebox(0,0)[cc]{$_{1_c}$}}
\put(90.33,17.00){\makebox(0,0)[cc]{$_{1_c}$}}
\put(116.33,46.00){\makebox(0,0)[cc]{$_{5_b}$}}
\put(104.33,59.00){\makebox(0,0)[cc]{$_{5_b}$}}
\put(101.33,52.00){\makebox(0,0)[cc]{$_{2_b}$}}
\put(96.33,59.00){\makebox(0,0)[cc]{$_{2_b}$}}
\put(60.33,46.00){\makebox(0,0)[cc]{$_{4_a}$}}
\put(69.33,53.00){\makebox(0,0)[cc]{$_{4_a}$}}
\put(80.33,61.00){\makebox(0,0)[cc]{$_{3_a}$}}
\put(94.33,46.00){\makebox(0,0)[cc]{$_{3_a}$}}
\put(87.33,93.00){\makebox(0,0)[cc]{$_{1_c}$}}
\put(86.33,72.00){\makebox(0,0)[cc]{$_{1_c}$}}
\put(83.33,81.00){\makebox(0,0)[cc]{$_{6_c}$}}
\put(90.33,81.00){\makebox(0,0)[cc]{$_{6_c}$}}
\put(112.33,63.00){\makebox(0,0)[cc]{$_{3_a}$}}
\put(102.33,44.00){\makebox(0,0)[cc]{$_{6_c}$}}
\put(72.33,43.00){\makebox(0,0)[cc]{$_{6_c}$}}
\put(90.33,32.00){\makebox(0,0)[cc]{$_{1_c}$}}
\put(82.33,51.00){\makebox(0,0)[cc]{$_{1_c}$}}
\put(74.33,68.00){\makebox(0,0)[cc]{$_{5_b}$}}
\put(99.33,77.00){\makebox(0,0)[cc]{$_{4_a}$}}
\put(96.33,67.00){\makebox(0,0)[cc]{$_{4_a}$}}
\put(97.33,87.00){\makebox(0,0)[cc]{$_{3_a}$}}
\put(78.33,87.00){\makebox(0,0)[cc]{$_{2_b}$}}
\put(63.33,91.00){\makebox(0,0)[cc]{$_{7_b}$}}
\put(51.33,70.00){\makebox(0,0)[cc]{$_{7_b}$}}
\put(111.33,90.00){\makebox(0,0)[cc]{$_{7_a}$}}
\put(123.33,74.00){\makebox(0,0)[cc]{$_{7_a}$}}
\put(113.33,78.00){\makebox(0,0)[cc]{$_{7_a}$}}
\put(77.33,28.00){\makebox(0,0)[cc]{$_{7_c}$}}
\put(77.33,34.00){\makebox(0,0)[cc]{$_{7_c}$}}
\put(100.33,27.00){\makebox(0,0)[cc]{$_{7_c}$}}
\put(67.33,89.00){\makebox(0,0)[cc]{$_{7_b}$}}
\put(56.33,70.00){\makebox(0,0)[cc]{$_{7_b}$}}
\put(99.33,33.00){\makebox(0,0)[cc]{$_{7_c}$}}
\put(119.33,68.00){\makebox(0,0)[cc]{$_{7_a}$}}
\put(107.33,87.00){\makebox(0,0)[cc]{$_{7_a}$}}
\put(73.33,80.00){\makebox(0,0)[cc]{$_{5_b}$}}
\put(132.33,1.00){\circle{2.00}} \put(178.33,55.00){\circle{2.00}}
\emline{114.33}{31.00}{93}{131.33}{2.00}{94}
\put(114.33,32.00){\circle{2.00}} \put(148.33,33.00){\circle{2.00}}
\emline{115.33}{32.00}{95}{147.33}{33.00}{96}
\emline{149.33}{34.00}{97}{177.33}{54.00}{98}
\emline{147.33}{32.00}{99}{132.33}{2.00}{100}
\emline{177.33}{55.00}{101}{115.33}{33.00}{102}
\emline{133.33}{2.00}{103}{178.33}{54.00}{104}
\put(182.33,53.00){\makebox(0,0)[cc]{$_{3^c}$}}
\put(129.33,1.00){\makebox(0,0)[cc]{$_{1^c}$}}
\put(112.33,30.00){\makebox(0,0)[cc]{$_{5^a}$}}
\put(133.33,30.00){\makebox(0,0)[cc]{$_{2_b}$}}
\put(158.33,27.00){\makebox(0,0)[cc]{$_{2_b}$}}
\put(122.33,14.00){\makebox(0,0)[cc]{$_{4_a}$}}
\put(152.33,49.00){\makebox(0,0)[cc]{$_{6_c}$}}
\put(134.33,11.00){\makebox(0,0)[cc]{$_{6_c}$}}
\put(161.33,45.00){\makebox(0,0)[cc]{$_{4_a}$}}
\put(151.33,30.00){\makebox(0,0)[cc]{$_{7^f}$}}
\put(60.33,32.00){\circle{2.00}} \put(4.33,52.00){\circle{2.00}}
\emline{43.33}{2.00}{105}{60.33}{31.00}{106}
\emline{5.33}{52.00}{107}{59.33}{32.00}{108}
\put(42.33,1.00){\circle{2.00}}
\emline{5.33}{51.00}{109}{22.33}{34.00}{110}
\put(23.33,33.00){\circle{2.00}}
\emline{24.33}{33.00}{111}{59.33}{32.00}{112}
\emline{24.33}{32.00}{113}{42.33}{2.00}{114}
\emline{4.33}{51.00}{115}{41.33}{2.00}{116}
\put(63.33,30.00){\makebox(0,0)[cc]{$_{4^e}$}}
\put(46.33,1.00){\makebox(0,0)[cc]{$_{1^a}$}}
\put(2.33,50.00){\makebox(0,0)[cc]{$_{2^b}$}}
\put(27.33,30.00){\makebox(0,0)[cc]{$_{7^f}$}}
\put(53.33,14.00){\makebox(0,0)[cc]{$_{5_b}$}}
\put(31.33,13.00){\makebox(0,0)[cc]{$_{3_a}$}}
\put(20.33,49.00){\makebox(0,0)[cc]{$_{6_c}$}}
\put(49.33,30.00){\makebox(0,0)[cc]{$_{3_a}$}}
\put(39.33,12.00){\makebox(0,0)[cc]{$_{6_c}$}}
\put(10.33,48.00){\makebox(0,0)[cc]{$_{5_b}$}}
\put(189.33,71.00){\circle{2.00}} \put(224.33,71.00){\circle{2.00}}
\emline{190.33}{71.00}{117}{223.33}{71.00}{118}
\put(213.33,1.00){\circle{2.00}} \put(205.33,40.00){\circle{2.00}}
\emline{190.33}{70.00}{119}{205.33}{41.00}{120}
\emline{206.33}{41.00}{121}{223.33}{70.00}{122}
\emline{205.33}{39.00}{123}{213.33}{2.00}{124}
\emline{189.33}{70.00}{125}{212.33}{2.00}{126}
\emline{213.33}{2.00}{127}{224.33}{70.00}{128}
\put(194.33,69.00){\makebox(0,0)[cc]{$_{5^f}$}}
\put(210.33,1.00){\makebox(0,0)[cc]{$_{6^f}$}}
\put(220.33,69.00){\makebox(0,0)[cc]{$_{4^f}$}}
\put(209.33,37.00){\makebox(0,0)[cc]{$_{7^f}$}}
\put(220.33,27.00){\makebox(0,0)[cc]{$_{2_b}$}}
\put(198.33,36.00){\makebox(0,0)[cc]{$_{3_a}$}}
\put(206.33,68.00){\makebox(0,0)[cc]{$_{1_c}$}}
\put(211.33,25.00){\makebox(0,0)[cc]{$_{1_c}$}}
\put(218.33,52.00){\makebox(0,0)[cc]{$_{4_a}$}}
\put(216.33,62.00){\makebox(0,0)[cc]{$_{3_a}$}}
\put(197.33,62.00){\makebox(0,0)[cc]{$_{2_b}$}}
\put(206.33,93.00){\circle{2.00}} \put(174.33,61.00){\circle{2.00}}
\put(147.33,98.00){\circle{2.00}} \put(179.33,86.00){\circle{2.00}}
\emline{148.33}{98.00}{129}{178.33}{86.00}{130}
\emline{180.33}{86.00}{131}{205.33}{92.00}{132}
\emline{179.33}{85.00}{133}{174.33}{62.00}{134}
\emline{148.33}{98.00}{135}{205.33}{93.00}{136}
\emline{206.33}{92.00}{137}{175.33}{62.00}{138}
\emline{173.33}{62.00}{139}{148.33}{97.00}{140}
\put(148.33,102.00){\makebox(0,0)[cc]{$_{2^d}$}}
\put(177.33,58.00){\makebox(0,0)[cc]{$_{6^d}$}}
\put(210.33,93.00){\makebox(0,0)[cc]{$_{3^d}$}}
\put(183.33,83.00){\makebox(0,0)[cc]{$_{7^f}$}}
\put(197.33,80.00){\makebox(0,0)[cc]{$_{5_b}$}}
\put(162.33,74.00){\makebox(0,0)[cc]{$_{4_a}$}}
\put(179.33,93.00){\makebox(0,0)[cc]{$_{1_c}$}}
\put(175.33,72.00){\makebox(0,0)[cc]{$_{1_c}$}}
\put(167.33,89.00){\makebox(0,0)[cc]{$_{5_b}$}}
\put(197.33,88.00){\makebox(0,0)[cc]{$_{4_a}$}}
\emline{139.00}{96.00}{141}{117.00}{87.00}{142}
\emline{112.00}{72.00}{143}{87.00}{65.00}{144}
\emline{53.00}{64.00}{145}{115.00}{86.00}{146}
\put(93.33,42.00){\makebox(0,0)[cc]{$_{1_c}$}}
\end{picture}
\caption{Disposition of copies of $K_{2,2,2}$ and $K_4$ at vertex
$7^f$ in $G$}
\end{figure}

Each point $p$ of ${\mathcal F}$ determines a {\it Pasch configuration}
$PC(p)$, formed by the 4 lines of ${\mathcal F}$ that
do not contain $p$. This $PC(p)$ may be denoted also
$pc(q_ar_a,q_br_b,q_cr_c),$ where $pq_ar_a,pq_br_b,pq_cr_c$ are the
lines of ${\mathcal F}$ containing $p$.
None of the lines of $PC(p)$
contains either $q_ar_a$ or $q_br_b$ or $q_cr_c$.
The 7 possible Pasch configurations here are:

$$\begin{array}{l}
PC(1)=pc(23,45,67)=\{246,257,347,356\}, \\
PC(2)=pc(13,46,57)=\{145,167,347,356\}, \\
PC(3)=pc(12,47,56)=\{145,167,246,257\}, \\
PC(4)=pc(15,26,37)=\{123,167,257,357\}, \\
PC(5)=pc(14,27,36)=\{123,167,246,347\}, \\
PC(6)=pc(17,24,35)=\{123,145,257,347\}, \\
PC(7)=pc(16,25,34)=\{123,145,246,356\}.
\end{array}$$
Figure 2 shows the disposition of the induced copies of $K_{2,2,2}$ and $K_4$
incident to the vertex $7^f$ in $G$,
represented by 3 octahedra and 4 tetrahedra, respectively, with
vertices and edges accompanied by their respective simplified notations and weak colors.
The 4 tetrahedra in the figure are also shown as separate entities,
for better distinction,
while the 3 octahedra are integrated in the central drawing as an upper-left,
an upper-right and a lower-central octahedron, radiated from the central vertex,
$7^f$. This 7 polyhedra can be blown up to 3-space without more
intersections than those of the vertices and edges shown in the figure.
Starting from the right upper corner in the figure and shown counterclockwise,
the 3 octahedra have respective composing 4-holes, each
sub-indexed with its common weak color, as follows:
$$\begin{array}{cccc}
{[347]}_a: & (7^f,4^e,7^e,4^f)_{3_a}, & (7^f,3^d,7^e,3^e)_{4_a}, & (3^d,4^e,3^e,4^f)_{7_a}; \\
{[257]}_b: & (7^f,5^a,7^a,5^f)_{2_b}, & (7^f,2^b,7^a,2^d)_{5_b}, & (2^b,5^a,2^d,5^f)_{7_b}; \\
{[167]}_c: & (7^f,6^d,7^d,6^f)_{1_c}, & (7^f,1^a,7^d,1^c)_{6_c}, & (1^a,6^d,1^c,6^f)_{7_c}.
\end{array}$$
The triangles in each octahedron here differ from those in the
copies of $K_4$ in $G$ in the way their edges are weakly colored. For example, the
copies of $K_4$ in Figure 2, namely those denoted $\langle 321\rangle $, $\langle 426\rangle $,
$\langle 356\rangle $, $\langle 451\rangle $, have their corresponding sets of constituent
triangles with the clockwise sequences of simplified notations and
weak colors of respective alternate incident vertices and edges, as
follows:
$$\begin{array}{cccc}
\{ (7^f,2_b,5^f,1_c,4^f,3_a), & (6^f,3_a,5^f,1_c,4^f,2_b), & %green
 (7^f,1_c,6^f,3_a,5^f,2_b), & (7^f,3_a,4^f,2_b,6^f,6_c)\}; \\
\{(7^f,6_c,1^c,4_a,5^a,2_b), & (3^e,2_b,1^c,4_a,5^a,2_b), & %red
 (7^f,4_a,3^e,2_b,1^c,4_a), & (7^f,2_b,5^a,6_c,3^e,4_a)\}; \\
\{(7^f,3_a,4^e,5_b,1^a,6_c), & (2^b,6_c,4^e,5_b,1^a,3_a), &
 (7^f,6_c,1^a,3_a,2^b,5_b), & (7^f,5_b,2^b,6_c,4^e,3_a)\}; \\ %blue
\{ (2^d,1_c,3^d,5_b,6^d,4_a), & (7^f,5_b,2^d,1_c,3^d,4_a), &
 (7^f,6_c,6^d,4_a,2^d,5_b), & (7^f,4_a,3^d,5_b,6^d,6_c)\}. %brown
\end{array}$$
This reflects the fact that the vertex $7^f=(7,34,25,16)=(p,q_a,r_a,q_br_b,q_cr_c)$
is associated with the Pasch configuration
$pc(34,25,16)=pc(q_ar_a,q_br_b,q_cr_c)$
given with its triples ordered according to the presence of the different
symbols $\ne 7$ at the 3 pair positions $a,b,c$, which is shown in
the ordered Fano lines $321=q_aq_bq_c$, $426=r_aq_br_c$,
$356=q_ar_br_c$, $451=r_ar_bq_c$, or in their respectively associated
tetrahedra: $\langle 321\rangle ,\langle 426\rangle ,\langle 356\rangle ,\langle 451\rangle $. These ordered lines form
the {\it ordered Pasch configuration} $\overline{pc}(7^f)=\{321,426,356,451\}$.
Similarly, an ordered Pasch configuration is associated to
the set of copies of $K_4$ incident to any other vertex of $G$.
Moreover, the following two results are readily checked.

\begin{teo} Any vertex $v=(p,q_ar_a,q_br_b,q_cr_c)$ of $G$ can be expressed in such a way that
$\langle q_aq_bq_c\rangle ,$ $\langle q_ar_br_c\rangle ,$ $\langle r_aq_br_c\rangle ,$ $\langle r_ar_bq_c\rangle $ are its 4 incident copies of
$K_4$, reflecting their notation and that of its 3 incident octahedra.\end{teo}

\begin{demo} The ordered Pasch configuration $\overline{pc}(v)$
associated to $v$ determines the ordered lines
$q_aq_bq_c,$ $q_ar_br_c,$ $r_aq_br_c,$ $r_ar_bq_c$ associated to the copies of $K_4$,
while the 3 remaining triples of ${\mathcal F}$ provide the data for the
octahedra incident to $v$: $[pq_ar_a]_a, [pq_br_b]_b, [pq_cr_c]_c$.
\end{demo}

\begin{teo} For any edge $e$ of $G$, there exists exactly one copy of $K_{2,2,2}$
and one of $K_4$ in $G$ that intersect at $e$. Moreover, $e$ is the only edge at which those
copies intersect. Thus, $G$ is fastened.
\end{teo}

\begin{demo} Let $e=vv'$ have weak color $q_j$, where $v=(p,q_ar_a,q_br_b,q_cr_c)$ and
$v'=(p',q'_ar'_a,q'_br'_b,$ $q'_cr'_c)$.
Then, the octahedron $[pp'p'']_j$ and the tetrahedron $\langle xyz\rangle $ are the copies of $K_{2,2,2}$
and $K_4$ in the statement,
where: {\bf(a)} $pp'p''$ is the Fano line containing $p$ and $p'$,
{\bf(b)} $j\in\{a,b,c\}$ is such that
$pp''=q'_jr'_j$ and $pp'=q_jr_j$ and
{\bf(c)} $xyz$,
one of the 4 ordered lines cited in Theorem 4.2 with respect to $v$,
is the strong color of $e$.\end{demo}

For example, the edge $7^f5^a$ has weak color $2_b$ and strong color
426. This is the only edge shared by the octahedron $[257]_b$ and the tetrahedron $\langle 426\rangle $.

\section{Symmetric properties of $G$}

Each auto\-mor\-phism $\tau\in{\mathcal A}(G)$, is the composition of a
permutation $\phi^\tau$ of ${\mathcal F}$ with a permutation $\psi^\tau$
of $\{a,b,c\}$. A set of 16 generators $\tau_i$ of ${\mathcal A}(G)$, ($i=1\ldots 16$),
is given by $\tau_i=\psi_i\circ\phi_i=\phi_i\circ\psi_i$, where we denote
$\phi_i=\phi^{\tau_i}$, $\psi_i=\psi^{\tau_i}$, and
$$\begin{array}{llll}
\phi_1=(23)(67), & \phi_2=(45)(67), &
\phi_3=(13)(57), & \phi_4=(46)(57), \\
\phi_5=(12)(56), & \phi_6=(47)(56), &
\phi_7=(15)(37), & \phi_8=(26)(37), \\
\phi_9=(14)(27), & \phi_{10}=(27)(36), &
\phi_{11}=(17)(35), & \phi_{12}=(24)(35), \\
\phi_{13}=(16)(34), & \phi_{14}=(25)(34), &
\psi_{15}=(ac), & \psi_{16}=(bc),
\end{array}$$
with $\psi_i$ and $\phi_j$ taken as the identity maps of ${\mathcal F}$ and
$\{a,b,c\}$, respectively, for $1\le i\le 14$ and $j=15,16$.

The subgroup $\Gamma\subset{\mathcal A}(G)$ that sends the
lexicographically smallest arc $(1^a,2^a)$ onto itself, either
directly or inversely oriented, includes exchanging, or not, its
incident triangles $(2^a,3^a,1^a)$ and $(2^a,3^b,1^a)$ in $[123]_a$,
or $(1^a,2^a,6^f)$ and $(1^a,2^a,5^e)$ in $\langle 347\rangle $. Thus, $\Gamma$
contains 4 elements and has generating set
$\{\tau_6\circ\tau_{16},\tau_5\}$. Moreover, $\Gamma$ is a subgroup
of ${\mathcal A}([1]_a)$, which has generating set
$\{\tau_1,\tau_2,\tau_5,\tau_6,\tau_{16}\}$. Furthermore,
$\{\tau_1,\tau_2,\tau_5,\tau_6,\tau_{15},\tau_{16}\}$ is a
generating set for ${\mathcal A}(\cup_{j=a}^{c}[123]_j)$. The remaining
automorphisms $\tau_i$ map ${\mathcal A}(\cup_{j=a}^{c}[123]_j)$ onto
its nontrivial cosets in ${\mathcal A}(G)$ by left multiplication. The
subgroup of ${\mathcal A}(G)$ that fixes $1^a$ has order 24 and
generating set
$\{\phi_1,\,\,\phi_2,\,\,\phi_4\circ\phi_{16},\,\,\phi_6\circ\phi_{16},\,\,
\phi_8\circ\phi_{15},\,\,\phi_{10}\circ\phi_{15},\,\,
\phi_{12}\circ\phi_{15}\circ\phi_{16},\,\,\phi_{14}\circ\phi_{15}\circ\phi_{16}\}.$

\begin{teo} $G$ is a fastened $\{K_4,K_{2,2,2}\}$-ultrahomogeneous
$\{K_4\}_{42}^4\{K_{2,2,2}\}_{21}^3$-graph which is non-line-graphical, with
$$|{\mathcal A}(G)|=1008=4|E(G)|.$$ The edges of $G$ can be seen as the left
cosets of a subgroup $\Gamma\subset{\mathcal A}(G)$ of order 4, and its
vertices as the left cosets of a subgroup of ${\mathcal A}(G)$ of order
24.
\end{teo}

\begin{demo}
Recall from Subsection 3.1 that $|{\mathcal A}(G)|=1008$.

Notice that ${\mathcal A}(\langle  347\rangle)=S_4$ is formed by 24 automorphisms. Since $\frac{|{\mathcal A}(G)|}{|{\mathcal A}(\langle 347\rangle )|}=\frac{1008}{24}=42$, then $\langle 347\rangle $ is sent by an automorphism of $G$ onto any other copy of $K_4$ in $G$, from which it is not difficult to see that $G$ is $K_4$-ultrahomogeneous.

Similarly, ${\mathcal A}([123]_a)$ is formed by 48 automorphisms. Since $\frac{|{\mathcal A}(G)|}{|{\mathcal A}([123]_a)|}=\frac{1008}{48}=21$, then $[123]_a$ is sent by an automorphism of $G$ onto any other copy of $K_{2,2,2}$ in $G$, from which it is not difficult to see that $G$ is $K_{2,2,2}$-ultrahomogeneous.

On the other hand, $|{\mathcal A}(K_{2,2,2})|=48$ and $|E(K_{2,2,2})|=12$
agree with the fact that $|\Gamma|=\frac{|{\mathcal
A}(K_{2,2,2})|}{|E(K_{2,2,2})|}=4$. Since $G$ is the edge-disjoint
union of 21 copies of $K_{2,2,2}$, it contains a total of
$21|E(K_{2,2,2}) |=21\times 12=252$ edges. Now, $|{\mathcal
A}(G)|=1008=21\times 48=21|{\mathcal A}(K_{2,2,2})|$.  This is 4 times the
number 252 of edges of $G$.
These edges correspond to the left cosets of
$\Gamma$ in ${\mathcal A}(G)$ and its vertices to the left cosets of the
stabilizer of $1^a=(1,23,45,67)$ in ${\mathcal A}(G)$,
whose order is 24.
\end{demo}

\section{Configurations associated with $G$}

The symmetrical disposition of objects in $G$ gives place to several combinatorial
point-line configurations and to their associated Levi, Menger, and dual Menger graphs.

We present the points and lines of 3 self-dual configurations
obtained from $G$, and their incidence relations:
\begin{enumerate}
\item the 42 vertices and 42 tetrahedra of $G$, and incidence given
by inclusion of a vertex in a tetrahedron; this is a self-dual
$(42_4)$-configuration with 2-arc-transitive Levi graph of diameter
$=$ girth $=6$, auto\-mor\-phism-group order 2016, stabilizer order
24, distance distribution vector $(1,4,12,24,27,14,2)$ and
isomorphic arc-transitive Menger graphs of diameter $=$ girth $=3$,
degree 12 and auto\-mor\-phism-group order 1008;

\item the 168 tetrahedral triangles and 168 octahedral triangles in $G$ and their
sharing of an edge;
this is a self-dual $(168_6)$-configuration with semisymmetric
Levi graph of diameter $=$ girth
$=6$, auto\-mor\-phism-group order 1008, common stabilizer order 6,
distance distribution vectors
$(1,6,24,60,111,102,32)$ and $(1,6,24,60,108,102,35)$, (just differing at distances 4 and 6
by 3 vertices)
and vertex-transitive Menger graphs of common
degree 24, diameter $=$ girth $=3$ and auto\-mor\-phism-group orders 1008 and 2016,
respectively.

\item the 168 tetrahedral triangles and 168 octahedral triangles in $G$ and their
sharing of an edge;
this is a self-dual $(168_6)$-configuration with semisymmetric
Levi graph of diameter $=$ girth
$=6$, auto\-mor\-phism-group order 1008, common stabilizer order 6,
distance distribution vectors
$(1,6,24,60,111,102,32)$ and $(1,6,24,60,108,102,35)$, (just differing at distances 4 and 6
by 3 vertices)
and vertex-transitive Menger graphs of common
degree 24, diameter $=$ girth $=3$ and auto\-mor\-phism-group orders 1008 and 2016,
respectively.
\end{enumerate}
Another interesting configuration associated to $G$ is formed by
the 42 tetrahedra and 21 octahedra of $G$, and their sharing of an edge; this
is a flag-transitive $(42_6,21_{12}$)-configuration.

\example Let $\mathcal L$ be the Levi graph of the
$(42_4)$-configuration in item 1 above. Then
$$((1,23,45,67), \langle 246\rangle , (3,12,47,56), \langle 145\rangle , (6,17,24,35), \langle 725\rangle , (1,67,23,45))$$
$$\mbox{and }(\langle 123\rangle ,(4,15,26,37),\langle 167\rangle ,(2,13,46,57),\langle 347\rangle ,(5,36,15,27),\langle 312\rangle )$$
are the lexicographically smallest paths realizing the diameter of
$\mathcal L$ and departing from each one of the two vertex parts
of $\mathcal L$. The second lexicographically smallest paths are
$$((1,23,45,67), \langle 246\rangle , (3,12,47,56), \langle 176\rangle , (4,15,37,26), \langle 572\rangle , (1,45,67,23))$$
$$\mbox{and }(\langle 123\rangle ,(4,15,26,37),\langle 167\rangle ,(3,12,56,47),\langle 264\rangle ,(5,27,36,14),\langle 231\rangle ).$$
We reach this way to the only two vertices realizing
the diameter of $\mathcal L$ starting from (1, 23, 45, 67), namely
$(1,67,23,45)$ and $(1,45,67,23)$; respectively: starting from $\langle 123\rangle $,
namely $\langle 312\rangle $ and $\langle 231\rangle $.
Those two pairs of paths reflect the correspondence between both
parts of $\mathcal L$ induced by the map $\Phi$ in Section 2.

\section{On 6-holes and other subgraphs of $G$}

\begin{figure}[htp]
\unitlength=0.70mm
\special{em:linewidth 0.4pt}
\linethickness{0.4pt}
\begin{picture}(193.00,66.00)
\put(45.00,24.00){\makebox(0,0)[cc]{$_1$}}
\put(55.00,24.00){\makebox(0,0)[cc]{$_2$}}
\put(65.00,24.00){\makebox(0,0)[cc]{$_3$}}
\put(56.00,17.00){\makebox(0,0)[cc]{$_{5_b}$}}
\put(49.00,15.00){\makebox(0,0)[cc]{$_6$}}
\put(61.00,15.00){\makebox(0,0)[cc]{$_4$}}
\put(55.00,6.00){\makebox(0,0)[cc]{$_7$}}
\emline{47.00}{24.00}{1}{53.00}{24.00}{2}
\emline{57.00}{24.00}{3}{63.00}{24.00}{4}
\emline{55.00}{15.00}{5}{55.00}{8.00}{6}
\emline{45.00}{22.00}{7}{48.00}{17.00}{8}
\emline{50.00}{13.00}{9}{53.00}{8.00}{10}
\emline{57.00}{19.00}{11}{63.00}{23.00}{12}
\emline{57.00}{8.00}{13}{60.00}{13.00}{14}
\emline{62.00}{17.00}{15}{65.00}{22.00}{16}
\emline{53.00}{19.00}{17}{47.00}{23.00}{18}
\put(55.00,17.00){\circle{10.00}}
\emline{54.00}{17.00}{19}{50.00}{15.00}{20}
\emline{56.00}{17.00}{21}{60.00}{15.00}{22}
\put(85.00,24.00){\makebox(0,0)[cc]{$_5$}}
\put(95.00,24.00){\makebox(0,0)[cc]{$_2$}}
\put(105.00,24.00){\makebox(0,0)[cc]{$_7$}}
\put(96.00,17.00){\makebox(0,0)[cc]{$_{1_b}$}}
\put(89.00,15.00){\makebox(0,0)[cc]{$_6$}}
\put(101.00,15.00){\makebox(0,0)[cc]{$_4$}}
\put(95.00,6.00){\makebox(0,0)[cc]{$_3$}}
\emline{87.00}{24.00}{23}{93.00}{24.00}{24}
\emline{97.00}{24.00}{25}{103.00}{24.00}{26}
\emline{95.00}{15.00}{27}{95.00}{8.00}{28}
\emline{85.00}{22.00}{29}{88.00}{17.00}{30}
\emline{90.00}{13.00}{31}{93.00}{8.00}{32}
\emline{97.00}{19.00}{33}{103.00}{23.00}{34}
\emline{97.00}{8.00}{35}{100.00}{13.00}{36}
\emline{102.00}{17.00}{37}{105.00}{22.00}{38}
\emline{93.00}{19.00}{39}{87.00}{23.00}{40}
\put(95.00,17.00){\circle{10.00}}
\emline{94.00}{17.00}{41}{90.00}{15.00}{42}
\emline{96.00}{17.00}{43}{100.00}{15.00}{44}
\put(125.00,24.00){\makebox(0,0)[cc]{$_7$}}
\put(135.00,24.00){\makebox(0,0)[cc]{$_2$}}
\put(145.00,24.00){\makebox(0,0)[cc]{$_5$}}
\put(136.00,17.00){\makebox(0,0)[cc]{$_{3_b}$}}
\put(129.00,15.00){\makebox(0,0)[cc]{$_6$}}
\put(141.00,15.00){\makebox(0,0)[cc]{$_4$}}
\put(135.00,6.00){\makebox(0,0)[cc]{$_1$}}
\emline{127.00}{24.00}{45}{133.00}{24.00}{46}
\emline{137.00}{24.00}{47}{143.00}{24.00}{48}
\emline{135.00}{15.00}{49}{135.00}{8.00}{50}
\emline{125.00}{22.00}{51}{128.00}{17.00}{52}
\emline{130.00}{13.00}{53}{133.00}{8.00}{54}
\emline{137.00}{19.00}{55}{143.00}{23.00}{56}
\emline{137.00}{8.00}{57}{140.00}{13.00}{58}
\emline{142.00}{17.00}{59}{145.00}{22.00}{60}
\emline{133.00}{19.00}{61}{127.00}{23.00}{62}
\put(135.00,17.00){\circle{10.00}}
\emline{134.00}{17.00}{63}{130.00}{15.00}{64}
\emline{136.00}{17.00}{65}{140.00}{15.00}{66}
\put(165.00,24.00){\makebox(0,0)[cc]{$_3$}}
\put(175.00,24.00){\makebox(0,0)[cc]{$_2$}}
\put(185.00,24.00){\makebox(0,0)[cc]{$_1$}}
\put(176.00,17.00){\makebox(0,0)[cc]{$_{7_b}$}}
\put(169.00,15.00){\makebox(0,0)[cc]{$_6$}}
\put(181.00,15.00){\makebox(0,0)[cc]{$_4$}}
\put(175.00,6.00){\makebox(0,0)[cc]{$_5$}}
\emline{167.00}{24.00}{67}{173.00}{24.00}{68}
\emline{177.00}{24.00}{69}{183.00}{24.00}{70}
\emline{175.00}{15.00}{71}{175.00}{8.00}{72}
\emline{165.00}{22.00}{73}{168.00}{17.00}{74}
\emline{170.00}{13.00}{75}{173.00}{8.00}{76}
\emline{177.00}{19.00}{77}{183.00}{23.00}{78}
\emline{177.00}{8.00}{79}{180.00}{13.00}{80}
\emline{182.00}{17.00}{81}{185.00}{22.00}{82}
\emline{173.00}{19.00}{83}{167.00}{23.00}{84}
\put(175.00,17.00){\circle{10.00}}
\emline{174.00}{17.00}{85}{170.00}{15.00}{86}
\emline{176.00}{17.00}{87}{180.00}{15.00}{88}
\emline{41.00}{25.00}{89}{54.00}{2.00}{90}
\put(41.00,26.00){\circle{2.00}}
\put(55.00,1.00){\circle{2.00}}
\emline{56.00}{2.00}{91}{69.00}{25.00}{92}
\put(69.00,26.00){\circle{2.00}}
\emline{42.00}{26.00}{93}{68.00}{26.00}{94}
\emline{81.00}{25.00}{95}{94.00}{2.00}{96}
\put(81.00,26.00){\circle{2.00}}
\put(95.00,1.00){\circle{2.00}}
\emline{96.00}{2.00}{97}{109.00}{25.00}{98}
\put(109.00,26.00){\circle{2.00}}
\emline{82.00}{26.00}{99}{108.00}{26.00}{100}
\emline{121.00}{25.00}{101}{134.00}{2.00}{102}
\put(121.00,26.00){\circle{2.00}}
\put(135.00,1.00){\circle{2.00}}
\emline{136.00}{2.00}{103}{149.00}{25.00}{104}
\put(149.00,26.00){\circle{2.00}}
\emline{122.00}{26.00}{105}{148.00}{26.00}{106}
\emline{161.00}{25.00}{107}{174.00}{2.00}{108}
\put(161.00,26.00){\circle{2.00}}
\put(175.00,1.00){\circle{2.00}}
\emline{176.00}{2.00}{109}{189.00}{25.00}{110}
\put(189.00,26.00){\circle{2.00}}
\emline{162.00}{26.00}{111}{188.00}{26.00}{112}
\emline{55.00}{22.00}{113}{55.00}{19.00}{114}
\emline{95.00}{22.00}{115}{95.00}{19.00}{116}
\emline{135.00}{22.00}{117}{135.00}{19.00}{118}
\emline{175.00}{22.00}{119}{175.00}{19.00}{120}
\put(73.00,26.00){\makebox(0,0)[cc]{$_{6^d}$}}
\put(78.00,26.00){\makebox(0,0)[cc]{$_{4^c}$}}
\put(55.00,29.00){\makebox(0,0)[cc]{$_{2_a}$}}
\put(45.00,12.00){\makebox(0,0)[cc]{$_{6_a}$}}
\put(66.00,12.00){\makebox(0,0)[cc]{$_{4_a}$}}
\put(59.00,1.00){\makebox(0,0)[cc]{$_{2^d}$}}
\put(113.00,26.00){\makebox(0,0)[cc]{$_{6^c}$}}
\put(95.00,29.00){\makebox(0,0)[cc]{$_{2_a}$}}
\put(85.00,12.00){\makebox(0,0)[cc]{$_{6_a}$}}
\put(106.00,12.00){\makebox(0,0)[cc]{$_{4_a}$}}
\put(99.00,1.00){\makebox(0,0)[cc]{$_{2^c}$}}
\put(38.00,26.00){\makebox(0,0)[cc]{$_{4^c}$}}
\put(153.00,26.00){\makebox(0,0)[cc]{$_{6^d}$}}
\put(158.00,26.00){\makebox(0,0)[cc]{$_{4^d}$}}
\put(135.00,29.00){\makebox(0,0)[cc]{$_{2_a}$}}
\put(125.00,12.00){\makebox(0,0)[cc]{$_{6_a}$}}
\put(146.00,12.00){\makebox(0,0)[cc]{$_{4_a}$}}
\put(139.00,1.00){\makebox(0,0)[cc]{$_{2^c}$}}
\put(193.00,26.00){\makebox(0,0)[cc]{$_{6^c}$}}
\put(175.00,29.00){\makebox(0,0)[cc]{$_{2_a}$}}
\put(165.00,12.00){\makebox(0,0)[cc]{$_{6_a}$}}
\put(186.00,12.00){\makebox(0,0)[cc]{$_{4_a}$}}
\put(179.00,1.00){\makebox(0,0)[cc]{$_{2^d}$}}
\put(118.00,26.00){\makebox(0,0)[cc]{$_{4^d}$}}
\put(45.00,61.00){\makebox(0,0)[cc]{$_1$}}
\put(55.00,61.00){\makebox(0,0)[cc]{$_2$}}
\put(65.00,61.00){\makebox(0,0)[cc]{$_3$}}
\put(56.00,54.00){\makebox(0,0)[cc]{$_{5_c}$}}
\put(49.00,52.00){\makebox(0,0)[cc]{$_6$}}
\put(61.00,52.00){\makebox(0,0)[cc]{$_4$}}
\put(55.00,43.00){\makebox(0,0)[cc]{$_7$}}
\emline{47.00}{61.00}{121}{53.00}{61.00}{122}
\emline{57.00}{61.00}{123}{63.00}{61.00}{124}
\emline{55.00}{52.00}{125}{55.00}{45.00}{126}
\emline{45.00}{59.00}{127}{48.00}{54.00}{128}
\emline{50.00}{50.00}{129}{53.00}{45.00}{130}
\emline{57.00}{56.00}{131}{63.00}{60.00}{132}
\emline{57.00}{45.00}{133}{60.00}{50.00}{134}
\emline{62.00}{54.00}{135}{65.00}{59.00}{136}
\emline{53.00}{56.00}{137}{47.00}{60.00}{138}
\put(55.00,54.00){\circle{10.00}}
\emline{54.00}{54.00}{139}{50.00}{52.00}{140}
\emline{56.00}{54.00}{141}{60.00}{52.00}{142}
\put(85.00,61.00){\makebox(0,0)[cc]{$_5$}}
\put(95.00,61.00){\makebox(0,0)[cc]{$_2$}}
\put(105.00,61.00){\makebox(0,0)[cc]{$_7$}}
\put(96.00,54.00){\makebox(0,0)[cc]{$_{1_c}$}}
\put(89.00,52.00){\makebox(0,0)[cc]{$_6$}}
\put(101.00,52.00){\makebox(0,0)[cc]{$_4$}}
\put(95.00,43.00){\makebox(0,0)[cc]{$_3$}}
\emline{87.00}{61.00}{143}{93.00}{61.00}{144}
\emline{97.00}{61.00}{145}{103.00}{61.00}{146}
\emline{95.00}{52.00}{147}{95.00}{45.00}{148}
\emline{85.00}{59.00}{149}{88.00}{54.00}{150}
\emline{90.00}{50.00}{151}{93.00}{45.00}{152}
\emline{97.00}{56.00}{153}{103.00}{60.00}{154}
\emline{97.00}{45.00}{155}{100.00}{50.00}{156}
\emline{102.00}{54.00}{157}{105.00}{59.00}{158}
\emline{93.00}{56.00}{159}{87.00}{60.00}{160}
\put(95.00,54.00){\circle{10.00}}
\emline{94.00}{54.00}{161}{90.00}{52.00}{162}
\emline{96.00}{54.00}{163}{100.00}{52.00}{164}
\put(125.00,61.00){\makebox(0,0)[cc]{$_7$}}
\put(135.00,61.00){\makebox(0,0)[cc]{$_2$}}
\put(145.00,61.00){\makebox(0,0)[cc]{$_5$}}
\put(136.00,54.00){\makebox(0,0)[cc]{$_{3_c}$}}
\put(129.00,52.00){\makebox(0,0)[cc]{$_6$}}
\put(141.00,52.00){\makebox(0,0)[cc]{$_4$}}
\put(135.00,43.00){\makebox(0,0)[cc]{$_1$}}
\emline{127.00}{61.00}{165}{133.00}{61.00}{166}
\emline{137.00}{61.00}{167}{143.00}{61.00}{168}
\emline{135.00}{52.00}{169}{135.00}{45.00}{170}
\emline{125.00}{59.00}{171}{128.00}{54.00}{172}
\emline{130.00}{50.00}{173}{133.00}{45.00}{174}
\emline{137.00}{56.00}{175}{143.00}{60.00}{176}
\emline{137.00}{45.00}{177}{140.00}{50.00}{178}
\emline{142.00}{54.00}{179}{145.00}{59.00}{180}
\emline{133.00}{56.00}{181}{127.00}{60.00}{182}
\put(135.00,54.00){\circle{10.00}}
\emline{134.00}{54.00}{183}{130.00}{52.00}{184}
\emline{136.00}{54.00}{185}{140.00}{52.00}{186}
\put(165.00,61.00){\makebox(0,0)[cc]{$_3$}}
\put(175.00,61.00){\makebox(0,0)[cc]{$_2$}}
\put(185.00,61.00){\makebox(0,0)[cc]{$_1$}}
\put(176.00,54.00){\makebox(0,0)[cc]{$_{7_c}$}}
\put(169.00,52.00){\makebox(0,0)[cc]{$_6$}}
\put(181.00,52.00){\makebox(0,0)[cc]{$_4$}}
\put(175.00,43.00){\makebox(0,0)[cc]{$_5$}}
\emline{167.00}{61.00}{187}{173.00}{61.00}{188}
\emline{177.00}{61.00}{189}{183.00}{61.00}{190}
\emline{175.00}{52.00}{191}{175.00}{45.00}{192}
\emline{165.00}{59.00}{193}{168.00}{54.00}{194}
\emline{170.00}{50.00}{195}{173.00}{45.00}{196}
\emline{177.00}{56.00}{197}{183.00}{60.00}{198}
\emline{177.00}{45.00}{199}{180.00}{50.00}{200}
\emline{182.00}{54.00}{201}{185.00}{59.00}{202}
\emline{173.00}{56.00}{203}{167.00}{60.00}{204}
\put(175.00,54.00){\circle{10.00}}
\emline{174.00}{54.00}{205}{170.00}{52.00}{206}
\emline{176.00}{54.00}{207}{180.00}{52.00}{208}
\emline{41.00}{62.00}{209}{54.00}{39.00}{210}
\put(41.00,63.00){\circle{2.00}}
\put(55.00,38.00){\circle{2.00}}
\emline{56.00}{39.00}{211}{69.00}{62.00}{212}
\put(69.00,63.00){\circle{2.00}}
\emline{42.00}{63.00}{213}{68.00}{63.00}{214}
\emline{81.00}{62.00}{215}{94.00}{39.00}{216}
\put(81.00,63.00){\circle{2.00}}
\put(95.00,38.00){\circle{2.00}}
\emline{96.00}{39.00}{217}{109.00}{62.00}{218}
\put(109.00,63.00){\circle{2.00}}
\emline{82.00}{63.00}{219}{108.00}{63.00}{220}
\emline{121.00}{62.00}{221}{134.00}{39.00}{222}
\put(121.00,63.00){\circle{2.00}}
\put(135.00,38.00){\circle{2.00}}
\emline{136.00}{39.00}{223}{149.00}{62.00}{224}
\put(149.00,63.00){\circle{2.00}}
\emline{122.00}{63.00}{225}{148.00}{63.00}{226}
\emline{161.00}{62.00}{227}{174.00}{39.00}{228}
\put(161.00,63.00){\circle{2.00}}
\put(175.00,38.00){\circle{2.00}}
\emline{176.00}{39.00}{229}{189.00}{62.00}{230}
\put(189.00,63.00){\circle{2.00}}
\emline{162.00}{63.00}{231}{188.00}{63.00}{232}
\emline{55.00}{59.00}{233}{55.00}{56.00}{234}
\emline{95.00}{59.00}{235}{95.00}{56.00}{236}
\emline{135.00}{59.00}{237}{135.00}{56.00}{238}
\emline{175.00}{59.00}{239}{175.00}{56.00}{240}
\put(73.00,63.00){\makebox(0,0)[cc]{$_{6^c}$}}
\put(78.00,63.00){\makebox(0,0)[cc]{$_{4^d}$}}
\put(55.00,66.00){\makebox(0,0)[cc]{$_{2_a}$}}
\put(45.00,49.00){\makebox(0,0)[cc]{$_{6_a}$}}
\put(66.00,49.00){\makebox(0,0)[cc]{$_{4_a}$}}
\put(59.00,38.00){\makebox(0,0)[cc]{$_{2^c}$}}
\put(113.00,63.00){\makebox(0,0)[cc]{$_{6^d}$}}
\put(95.00,66.00){\makebox(0,0)[cc]{$_{2_a}$}}
\put(85.00,49.00){\makebox(0,0)[cc]{$_{6_a}$}}
\put(106.00,49.00){\makebox(0,0)[cc]{$_{4_a}$}}
\put(99.00,38.00){\makebox(0,0)[cc]{$_{2^d}$}}
\put(38.00,63.00){\makebox(0,0)[cc]{$_{4^d}$}}
\put(153.00,63.00){\makebox(0,0)[cc]{$_{6^c}$}}
\put(158.00,63.00){\makebox(0,0)[cc]{$_{4^c}$}}
\put(135.00,66.00){\makebox(0,0)[cc]{$_{2_a}$}}
\put(125.00,49.00){\makebox(0,0)[cc]{$_{6_a}$}}
\put(146.00,49.00){\makebox(0,0)[cc]{$_{4_a}$}}
\put(139.00,38.00){\makebox(0,0)[cc]{$_{2^d}$}}
\put(193.00,63.00){\makebox(0,0)[cc]{$_{6^d}$}}
\put(175.00,66.00){\makebox(0,0)[cc]{$_{2_a}$}}
\put(165.00,49.00){\makebox(0,0)[cc]{$_{6_a}$}}
\put(186.00,49.00){\makebox(0,0)[cc]{$_{4_a}$}}
\put(179.00,38.00){\makebox(0,0)[cc]{$_{2^c}$}}
\put(118.00,63.00){\makebox(0,0)[cc]{$_{4^c}$}}
\end{picture}
\caption{Octahedral triangles of $[246]_a$}
\end{figure}
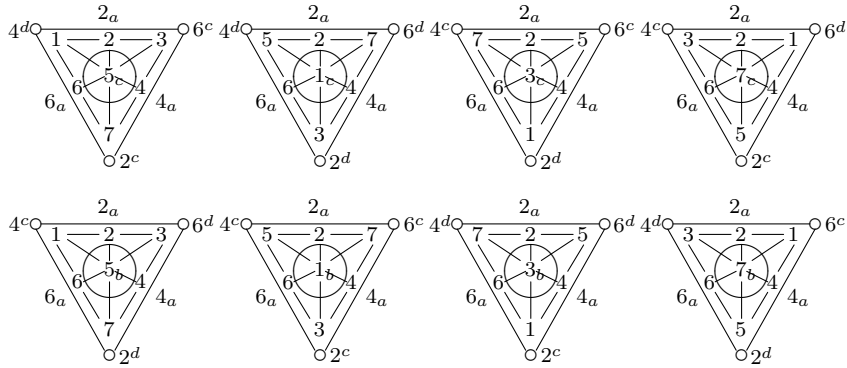

Let us depict the Fano plane in an octahedral triangle, as shown for example in Figure 3 for the 8 triangles of $[246]_a$, as follows.
For each such triangle $t$, there is a unique $i\in\{a,b,c\}$ and a unique point $p\in\mathcal F$ such that the 3 vertices of $t$ (looked upon as ordered pencils) have $p$ present in position $i$. The central point of the depiction of $\mathcal F$ inside each such $t$ is set to be $p$, subindexed by $i$.
For example in Figure 3, this $i$ appears as subindices $c$ and $b$, respectively, in the 4 top and 4 bottom triangles. Also, for each edge $e$ of $t$ with weak color $q_j$ and strong color $\ell$, the point in
the depiction of $\mathcal F$ in $t$ at the middle of $e$ is set to be $q$, and the 'external' line of $\mathcal F$ containing $q$ is set to be $\ell$, with the points of $\ell\setminus\{q\}$ set near the endvertices of $e$.

\begin{figure}[htp]
\unitlength=0.70mm
\special{em:linewidth 0.4pt}
\linethickness{0.4pt}
\begin{picture}(231.00,102.00)
\put(42.00,24.00){\makebox(0,0)[cc]{$_1$}}
\put(52.00,24.00){\makebox(0,0)[cc]{$_2$}}
\put(62.00,24.00){\makebox(0,0)[cc]{$_3$}}
\put(53.00,17.00){\makebox(0,0)[cc]{$_{5_c}$}}
\put(46.00,15.00){\makebox(0,0)[cc]{$_6$}}
\put(58.00,15.00){\makebox(0,0)[cc]{$_4$}}
\put(52.00,6.00){\makebox(0,0)[cc]{$_7$}}
\emline{44.00}{24.00}{1}{50.00}{24.00}{2}
\emline{54.00}{24.00}{3}{60.00}{24.00}{4}
\emline{52.00}{15.00}{5}{52.00}{8.00}{6}
\emline{42.00}{22.00}{7}{45.00}{17.00}{8}
\emline{47.00}{13.00}{9}{50.00}{8.00}{10}
\emline{54.00}{19.00}{11}{60.00}{23.00}{12}
\emline{54.00}{8.00}{13}{57.00}{13.00}{14}
\emline{59.00}{17.00}{15}{62.00}{22.00}{16}
\emline{50.00}{19.00}{17}{44.00}{23.00}{18}
\put(52.00,17.00){\circle{10.00}}
\emline{51.00}{17.00}{19}{47.00}{15.00}{20}
\emline{53.00}{17.00}{21}{57.00}{15.00}{22}
\emline{38.00}{25.00}{23}{51.00}{2.00}{24}
\put(38.00,26.00){\circle{2.00}}
\put(52.00,1.00){\circle{2.00}}
\emline{53.00}{2.00}{25}{66.00}{25.00}{26}
\put(66.00,26.00){\circle{2.00}}
\emline{39.00}{26.00}{27}{65.00}{26.00}{28}
\emline{52.00}{22.00}{29}{52.00}{19.00}{30}
\emline{24.00}{50.00}{31}{37.00}{27.00}{32}
\put(24.00,51.00){\circle{2.00}}
\put(10.00,26.00){\circle{2.00}}
\emline{11.00}{26.00}{33}{37.00}{26.00}{34}
\emline{11.00}{27.00}{35}{24.00}{50.00}{36}
\emline{10.00}{75.00}{37}{23.00}{52.00}{38}
\put(10.00,76.00){\circle{2.00}}
\put(38.00,76.00){\circle{2.00}}
\emline{11.00}{76.00}{39}{37.00}{76.00}{40}
\put(66.00,76.00){\circle{2.00}}
\emline{52.00}{100.00}{41}{65.00}{77.00}{42}
\put(52.00,101.00){\circle{2.00}}
\emline{39.00}{76.00}{43}{65.00}{76.00}{44}
\emline{39.00}{77.00}{45}{52.00}{100.00}{46}
\put(70.00,74.00){\makebox(0,0)[cc]{$_1$}}
\put(80.00,74.00){\makebox(0,0)[cc]{$_7$}}
\put(90.00,74.00){\makebox(0,0)[cc]{$_6$}}
\put(81.00,67.00){\makebox(0,0)[cc]{$_{5_c}$}}
\put(74.00,65.00){\makebox(0,0)[cc]{$_3$}}
\put(86.00,65.00){\makebox(0,0)[cc]{$_4$}}
\put(80.00,56.00){\makebox(0,0)[cc]{$_2$}}
\emline{72.00}{74.00}{47}{78.00}{74.00}{48}
\emline{82.00}{74.00}{49}{88.00}{74.00}{50}
\emline{80.00}{65.00}{51}{80.00}{58.00}{52}
\emline{70.00}{72.00}{53}{73.00}{67.00}{54}
\emline{75.00}{63.00}{55}{78.00}{58.00}{56}
\emline{82.00}{69.00}{57}{88.00}{73.00}{58}
\emline{82.00}{58.00}{59}{85.00}{63.00}{60}
\emline{87.00}{67.00}{61}{90.00}{72.00}{62}
\emline{78.00}{69.00}{63}{72.00}{73.00}{64}
\put(80.00,67.00){\circle{10.00}}
\emline{79.00}{67.00}{65}{75.00}{65.00}{66}
\emline{81.00}{67.00}{67}{85.00}{65.00}{68}
\emline{66.00}{75.00}{69}{79.00}{52.00}{70}
\put(80.00,51.00){\circle{2.00}}
\emline{81.00}{52.00}{71}{94.00}{75.00}{72}
\put(94.00,76.00){\circle{2.00}}
\emline{67.00}{76.00}{73}{93.00}{76.00}{74}
\emline{80.00}{72.00}{75}{80.00}{69.00}{76}
\put(94.00,26.00){\circle{2.00}}
\emline{80.00}{50.00}{77}{93.00}{27.00}{78}
\emline{67.00}{26.00}{79}{93.00}{26.00}{80}
\emline{67.00}{27.00}{81}{80.00}{50.00}{82}
\put(14.00,74.00){\makebox(0,0)[cc]{$_4$}}
\put(24.00,74.00){\makebox(0,0)[cc]{$_7$}}
\put(34.00,74.00){\makebox(0,0)[cc]{$_3$}}
\put(25.00,67.00){\makebox(0,0)[cc]{$_{5_c}$}}
\put(18.00,65.00){\makebox(0,0)[cc]{$_6$}}
\put(30.00,65.00){\makebox(0,0)[cc]{$_1$}}
\put(24.00,56.00){\makebox(0,0)[cc]{$_2$}}
\emline{16.00}{74.00}{83}{22.00}{74.00}{84}
\emline{26.00}{74.00}{85}{32.00}{74.00}{86}
\emline{24.00}{65.00}{87}{24.00}{58.00}{88}
\emline{14.00}{72.00}{89}{17.00}{67.00}{90}
\emline{19.00}{63.00}{91}{22.00}{58.00}{92}
\emline{26.00}{69.00}{93}{32.00}{73.00}{94}
\emline{26.00}{58.00}{95}{29.00}{63.00}{96}
\emline{31.00}{67.00}{97}{34.00}{72.00}{98}
\emline{22.00}{69.00}{99}{16.00}{73.00}{100}
\put(24.00,67.00){\circle{10.00}}
\emline{23.00}{67.00}{101}{19.00}{65.00}{102}
\emline{25.00}{67.00}{103}{29.00}{65.00}{104}
\emline{25.00}{52.00}{105}{38.00}{75.00}{106}
\emline{24.00}{72.00}{107}{24.00}{69.00}{108}
\put(24.00,46.00){\makebox(0,0)[cc]{$_2$}}
\put(25.00,34.00){\makebox(0,0)[cc]{$_{5_c}$}}
\put(14.00,27.00){\makebox(0,0)[cc]{$_6$}}
\put(24.00,27.00){\makebox(0,0)[cc]{$_7$}}
\put(34.00,27.00){\makebox(0,0)[cc]{$_1$}}
\emline{16.00}{28.00}{109}{22.00}{28.00}{110}
\emline{26.00}{28.00}{111}{32.00}{28.00}{112}
\put(24.00,34.00){\circle{10.00}}
\put(18.00,37.00){\makebox(0,0)[cc]{$_4$}}
\put(30.00,37.00){\makebox(0,0)[cc]{$_3$}}
\emline{15.00}{30.00}{113}{18.00}{35.00}{114}
\emline{20.00}{39.00}{115}{23.00}{44.00}{116}
\emline{25.00}{44.00}{117}{28.00}{39.00}{118}
\emline{30.00}{35.00}{119}{33.00}{30.00}{120}
\emline{24.00}{44.00}{121}{24.00}{37.00}{122}
\emline{16.00}{29.00}{123}{22.00}{33.00}{124}
\emline{32.00}{29.00}{125}{26.00}{33.00}{126}
\emline{24.00}{32.00}{127}{24.00}{29.00}{128}
\emline{19.00}{36.00}{129}{23.00}{34.00}{130}
\emline{29.00}{36.00}{131}{25.00}{34.00}{132}
\put(52.00,96.00){\makebox(0,0)[cc]{$_7$}}
\put(53.00,84.00){\makebox(0,0)[cc]{$_{5_c}$}}
\put(42.00,77.00){\makebox(0,0)[cc]{$_3$}}
\put(52.00,77.00){\makebox(0,0)[cc]{$_2$}}
\put(62.00,77.00){\makebox(0,0)[cc]{$_1$}}
\emline{44.00}{78.00}{133}{50.00}{78.00}{134}
\emline{54.00}{78.00}{135}{60.00}{78.00}{136}
\put(52.00,84.00){\circle{10.00}}
\put(46.00,87.00){\makebox(0,0)[cc]{$_4$}}
\put(58.00,87.00){\makebox(0,0)[cc]{$_6$}}
\emline{43.00}{80.00}{137}{46.00}{85.00}{138}
\emline{48.00}{89.00}{139}{51.00}{94.00}{140}
\emline{53.00}{94.00}{141}{56.00}{89.00}{142}
\emline{58.00}{85.00}{143}{61.00}{80.00}{144}
\emline{52.00}{94.00}{145}{52.00}{87.00}{146}
\emline{44.00}{79.00}{147}{50.00}{83.00}{148}
\emline{60.00}{79.00}{149}{54.00}{83.00}{150}
\emline{52.00}{82.00}{151}{52.00}{79.00}{152}
\emline{47.00}{86.00}{153}{51.00}{84.00}{154}
\emline{57.00}{86.00}{155}{53.00}{84.00}{156}
\put(80.00,46.00){\makebox(0,0)[cc]{$_2$}}
\put(81.00,34.00){\makebox(0,0)[cc]{$_{5_c}$}}
\put(70.00,27.00){\makebox(0,0)[cc]{$_3$}}
\put(80.00,27.00){\makebox(0,0)[cc]{$_7$}}
\put(90.00,27.00){\makebox(0,0)[cc]{$_4$}}
\emline{72.00}{28.00}{157}{78.00}{28.00}{158}
\emline{82.00}{28.00}{159}{88.00}{28.00}{160}
\put(80.00,34.00){\circle{10.00}}
\put(74.00,37.00){\makebox(0,0)[cc]{$_1$}}
\put(86.00,37.00){\makebox(0,0)[cc]{$_6$}}
\emline{71.00}{30.00}{161}{74.00}{35.00}{162}
\emline{76.00}{39.00}{163}{79.00}{44.00}{164}
\emline{81.00}{44.00}{165}{84.00}{39.00}{166}
\emline{86.00}{35.00}{167}{89.00}{30.00}{168}
\emline{80.00}{44.00}{169}{80.00}{37.00}{170}
\emline{72.00}{29.00}{171}{78.00}{33.00}{172}
\emline{88.00}{29.00}{173}{82.00}{33.00}{174}
\emline{80.00}{32.00}{175}{80.00}{29.00}{176}
\emline{75.00}{36.00}{177}{79.00}{34.00}{178}
\emline{85.00}{36.00}{179}{81.00}{34.00}{180}
\put(52.00,29.00){\makebox(0,0)[cc]{$_{2_a}$}}
\put(70.00,38.00){\makebox(0,0)[cc]{$_{1_b}$}}
\put(70.00,62.00){\makebox(0,0)[cc]{$_{3_a}$}}
\put(35.00,62.00){\makebox(0,0)[cc]{$_{1_a}$}}
\put(52.00,73.00){\makebox(0,0)[cc]{$_{2_b}$}}
\put(35.00,38.00){\makebox(0,0)[cc]{$_{3_b}$}}
\put(64.00,29.00){\makebox(0,0)[cc]{$_{6^c}$}}
\put(64.00,73.00){\makebox(0,0)[cc]{$_{4^f}$}}
\put(41.00,73.00){\makebox(0,0)[cc]{$_{6^a}$}}
\put(41.00,29.00){\makebox(0,0)[cc]{$_{4^d}$}}
\put(28.00,51.00){\makebox(0,0)[cc]{$_{7^b}$}}
\put(77.00,51.00){\makebox(0,0)[cc]{$_{7_e}$}}
\put(24.00,78.00){\makebox(0,0)[cc]{$_{7_a}$}}
\put(80.00,78.00){\makebox(0,0)[cc]{$_{7_a}$}}
\put(24.00,23.00){\makebox(0,0)[cc]{$_{7_b}$}}
\put(80.00,23.00){\makebox(0,0)[cc]{$_{7_b}$}}
\put(14.00,38.00){\makebox(0,0)[cc]{$_{4_b}$}}
\put(14.00,62.00){\makebox(0,0)[cc]{$_{6_a}$}}
\put(42.00,12.00){\makebox(0,0)[cc]{$_{6_a}$}}
\put(63.00,12.00){\makebox(0,0)[cc]{$_{4_a}$}}
\put(91.00,62.00){\makebox(0,0)[cc]{$_{4_a}$}}
\put(62.00,88.00){\makebox(0,0)[cc]{$_{6_b}$}}
\put(91.00,38.00){\makebox(0,0)[cc]{$_{6_b}$}}
\put(42.00,88.00){\makebox(0,0)[cc]{$_{4_b}$}}
\put(56.00,101.00){\makebox(0,0)[cc]{$_{2^a}$}}
\put(56.00,1.00){\makebox(0,0)[cc]{$_{2^c}$}}
\put(7.00,76.00){\makebox(0,0)[cc]{$_{1^e}$}}
\put(7.00,26.00){\makebox(0,0)[cc]{$_{3^a}$}}
\put(98.00,76.00){\makebox(0,0)[cc]{$_{3^c}$}}
\put(98.00,26.00){\makebox(0,0)[cc]{$_{1^b}$}}
\put(170.00,49.00){\circle{2.00}}
\put(150.00,49.00){\circle{2.00}}
\put(140.00,65.00){\circle{2.00}}
\put(180.00,65.00){\circle{2.00}}
\put(170.00,81.00){\circle{2.00}}
\put(200.00,33.00){\circle{2.00}}
\put(210.00,49.00){\circle{2.00}}
\put(200.00,65.00){\circle{2.00}}
\emline{151.00}{81.00}{181}{169.00}{81.00}{182}
\emline{171.00}{80.00}{183}{179.00}{66.00}{184}
\emline{179.00}{64.00}{185}{171.00}{50.00}{186}
\emline{169.00}{49.00}{187}{151.00}{49.00}{188}
\emline{149.00}{50.00}{189}{141.00}{64.00}{190}
\emline{141.00}{66.00}{191}{149.00}{80.00}{192}
\emline{181.00}{65.00}{193}{199.00}{65.00}{194}
\emline{201.00}{64.00}{195}{209.00}{50.00}{196}
\emline{209.00}{48.00}{197}{201.00}{34.00}{198}
\emline{199.00}{33.00}{199}{181.00}{33.00}{200}
\emline{179.00}{34.00}{201}{171.00}{48.00}{202}
\put(170.00,17.00){\circle{2.00}}
\put(150.00,17.00){\circle{2.00}}
\put(140.00,33.00){\circle{2.00}}
\emline{179.00}{32.00}{203}{171.00}{18.00}{204}
\emline{169.00}{17.00}{205}{151.00}{17.00}{206}
\emline{149.00}{18.00}{207}{141.00}{32.00}{208}
\emline{141.00}{34.00}{209}{149.00}{48.00}{210}
\put(200.00,1.00){\circle{2.00}}
\put(180.00,1.00){\circle{2.00}}
\put(210.00,17.00){\circle{2.00}}
\emline{201.00}{32.00}{211}{209.00}{18.00}{212}
\emline{209.00}{16.00}{213}{201.00}{2.00}{214}
\emline{199.00}{1.00}{215}{181.00}{1.00}{216}
\emline{179.00}{2.00}{217}{171.00}{16.00}{218}
\put(140.00,1.00){\circle{2.00}}
\emline{149.00}{16.00}{219}{141.00}{2.00}{220}
\put(180.00,97.00){\circle{2.00}}
\put(200.00,97.00){\circle{2.00}}
\emline{179.00}{96.00}{221}{171.00}{82.00}{222}
\emline{181.00}{97.00}{223}{199.00}{97.00}{224}
\emline{201.00}{96.00}{225}{209.00}{82.00}{226}
\emline{209.00}{80.00}{227}{201.00}{66.00}{228}
\put(230.00,49.00){\circle{2.00}}
\emline{211.00}{49.00}{229}{229.00}{49.00}{230}
\put(230.00,81.00){\circle{2.00}}
\emline{211.00}{81.00}{231}{229.00}{81.00}{232}
\emline{111.00}{80.00}{233}{119.00}{66.00}{234}
\put(120.00,97.00){\circle{2.00}}
\emline{119.00}{96.00}{235}{111.00}{82.00}{236}
\emline{121.00}{97.00}{237}{139.00}{97.00}{238}
\put(140.00,97.00){\circle{2.00}}
\emline{141.00}{96.00}{239}{149.00}{82.00}{240}
\emline{205.00}{97.00}{241}{211.00}{97.00}{242}
\emline{217.00}{97.00}{243}{223.00}{97.00}{244}
\emline{145.00}{97.00}{245}{151.00}{97.00}{246}
\emline{157.00}{97.00}{247}{163.00}{97.00}{248}
\emline{169.00}{97.00}{249}{175.00}{97.00}{250}
\put(110.00,49.00){\circle{2.00}}
\put(120.00,65.00){\circle{2.00}}
\emline{119.00}{64.00}{251}{111.00}{50.00}{252}
\emline{121.00}{65.00}{253}{139.00}{65.00}{254}
\emline{139.00}{33.00}{255}{121.00}{33.00}{256}
\emline{119.00}{34.00}{257}{111.00}{48.00}{258}
\put(110.00,17.00){\circle{2.00}}
\emline{119.00}{32.00}{259}{111.00}{18.00}{260}
\put(120.00,1.00){\circle{2.00}}
\emline{139.00}{1.00}{261}{121.00}{1.00}{262}
\emline{119.00}{2.00}{263}{111.00}{16.00}{264}
\emline{229.00}{17.00}{265}{211.00}{17.00}{266}
\emline{157.00}{1.00}{267}{163.00}{1.00}{268}
\emline{219.00}{1.00}{269}{225.00}{1.00}{270}
\put(173.00,40.00){\makebox(0,0)[cc]{$_{4_b}$}}
\put(143.00,74.00){\makebox(0,0)[cc]{$_{6_c}$}}
\put(208.00,58.00){\makebox(0,0)[cc]{$_{4_a}$}}
\put(179.00,88.00){\makebox(0,0)[cc]{$_{6_b}$}}
\put(160.00,78.00){\makebox(0,0)[cc]{$_{2_a}$}}
\put(160.00,20.00){\makebox(0,0)[cc]{$_{2_b}$}}
\put(179.00,9.00){\makebox(0,0)[cc]{$_{4_a}$}}
\put(173.00,26.00){\makebox(0,0)[cc]{$_{6_c}$}}
\put(190.00,30.00){\makebox(0,0)[cc]{$_{2_a}$}}
\put(203.00,23.00){\makebox(0,0)[cc]{$_{4_c}$}}
\put(203.00,42.00){\makebox(0,0)[cc]{$_{6_b}$}}
\put(219.00,46.00){\makebox(0,0)[cc]{$_{2_c}$}}
\put(202.00,9.00){\makebox(0,0)[cc]{$_{6_a}$}}
\put(190.00,94.00){\makebox(0,0)[cc]{$_{2_c}$}}
\put(113.00,89.00){\makebox(0,0)[cc]{$_{6_b}$}}
\put(208.00,90.00){\makebox(0,0)[cc]{$_{4_b}$}}
\put(129.00,94.00){\makebox(0,0)[cc]{$_{2_c}$}}
\put(220.00,90.00){\makebox(0,0)[cc]{$_{246_c^5}$}}
\put(220.00,65.00){\makebox(0,0)[cc]{$_{246_b^7}$}}
\put(190.00,49.00){\makebox(0,0)[cc]{$_{246_c^3}$}}
\put(160.00,7.00){\makebox(0,0)[cc]{$_{246_c^5}$}}
\put(220.00,33.00){\makebox(0,0)[cc]{$_{246_a^5}$}}
\put(190.00,4.00){\makebox(0,0)[cc]{$_{2_c}$}}
\put(160.00,52.00){\makebox(0,0)[cc]{$_{2_c}$}}
\put(190.00,81.00){\makebox(0,0)[cc]{$_{246_a^1}$}}
\put(190.00,17.00){\makebox(0,0)[cc]{$_{246_b^1}$}}
\put(160.00,33.00){\makebox(0,0)[cc]{$_{246_a^7}$}}
\put(160.00,65.00){\makebox(0,0)[cc]{$_{246_b^5}$}}
\put(113.00,71.00){\makebox(0,0)[cc]{$_{4_c}$}}
\put(220.00,78.00){\makebox(0,0)[cc]{$_{2_a}$}}
\put(190.00,62.00){\makebox(0,0)[cc]{$_{2_b}$}}
\put(177.00,75.00){\makebox(0,0)[cc]{$_{4_c}$}}
\put(179.00,55.00){\makebox(0,0)[cc]{$_{6_a}$}}
\put(148.00,89.00){\makebox(0,0)[cc]{$_{4_b}$}}
\put(148.00,58.00){\makebox(0,0)[cc]{$_{4_a}$}}
\put(160.00,90.00){\makebox(0,0)[cc]{$_{246_c^7}$}}
\put(148.00,38.00){\makebox(0,0)[cc]{$_{6_b}$}}
\put(148.00,27.00){\makebox(0,0)[cc]{$_{4_c}$}}
\put(148.00,6.00){\makebox(0,0)[cc]{$_{6_a}$}}
\put(208.00,70.00){\makebox(0,0)[cc]{$_{6_c}$}}
\put(130.00,81.00){\makebox(0,0)[cc]{$_{246_a^3}$}}
\put(130.00,49.00){\makebox(0,0)[cc]{$_{246_c^1}$}}
\put(130.00,17.00){\makebox(0,0)[cc]{$_{246_b^3}$}}
\put(220.00,7.00){\makebox(0,0)[cc]{$_{246_c^7}$}}
\put(130.00,62.00){\makebox(0,0)[cc]{$_{2_b}$}}
\put(118.00,55.00){\makebox(0,0)[cc]{$_{6_a}$}}
\put(118.00,43.00){\makebox(0,0)[cc]{$_{4_b}$}}
\put(118.00,23.00){\makebox(0,0)[cc]{$_{6_c}$}}
\put(118.00,10.00){\makebox(0,0)[cc]{$_{4_a}$}}
\put(222.00,20.00){\makebox(0,0)[cc]{$_{2_b}$}}
\put(130.00,4.00){\makebox(0,0)[cc]{$_{2_c}$}}
\put(130.00,36.00){\makebox(0,0)[cc]{$_{2_a}$}}
\put(150.00,81.00){\circle*{2.00}}
\put(120.00,33.00){\circle*{2.00}}
\put(210.00,81.00){\circle*{2.00}}
\put(180.00,33.00){\circle*{2.00}}
\put(110.00,81.00){\circle{2.00}}
\emline{230.00}{52.00}{271}{230.00}{56.00}{272}
\emline{230.00}{59.00}{273}{230.00}{63.00}{274}
\emline{230.00}{66.00}{275}{230.00}{70.00}{276}
\emline{230.00}{73.00}{277}{230.00}{77.00}{278}
\emline{230.00}{82.00}{279}{230.00}{86.00}{280}
\emline{230.00}{89.00}{281}{230.00}{93.00}{282}
\emline{227.00}{97.00}{283}{230.00}{97.00}{284}
\emline{110.00}{52.00}{285}{110.00}{56.00}{286}
\emline{110.00}{59.00}{287}{110.00}{63.00}{288}
\emline{110.00}{66.00}{289}{110.00}{70.00}{290}
\emline{110.00}{73.00}{291}{110.00}{77.00}{292}
\emline{110.00}{82.00}{293}{110.00}{86.00}{294}
\emline{110.00}{89.00}{295}{110.00}{93.00}{296}
\emline{116.00}{97.00}{297}{119.00}{97.00}{298}
\emline{110.00}{97.00}{299}{113.00}{97.00}{300}
\put(230.00,17.00){\circle{2.00}}
\emline{230.00}{20.00}{301}{230.00}{24.00}{302}
\emline{230.00}{27.00}{303}{230.00}{31.00}{304}
\emline{230.00}{34.00}{305}{230.00}{38.00}{306}
\emline{230.00}{41.00}{307}{230.00}{45.00}{308}
\emline{230.00}{2.00}{309}{230.00}{6.00}{310}
\emline{230.00}{9.00}{311}{230.00}{13.00}{312}
\emline{110.00}{20.00}{313}{110.00}{24.00}{314}
\emline{110.00}{27.00}{315}{110.00}{31.00}{316}
\emline{110.00}{34.00}{317}{110.00}{38.00}{318}
\emline{110.00}{41.00}{319}{110.00}{45.00}{320}
\emline{110.00}{2.00}{321}{110.00}{6.00}{322}
\emline{110.00}{9.00}{323}{110.00}{13.00}{324}
\put(52.00,51.00){\makebox(0,0)[cc]{$_{123_c^5}$}}
\put(91.00,50.00){\makebox(0,0)[cc]{$_{246_c^5}$}}
\put(14.00,51.00){\makebox(0,0)[cc]{$_{246_c^5}$}}
\put(24.00,88.00){\makebox(0,0)[cc]{$_{347_c^5}$}}
\put(80.00,88.00){\makebox(0,0)[cc]{$_{167_c^5}$}}
\put(24.00,12.00){\makebox(0,0)[cc]{$_{167_c^5}$}}
\put(80.00,12.00){\makebox(0,0)[cc]{$_{347_c^5}$}}
\emline{230.00}{2.00}{325}{230.00}{1.00}{326}
\emline{110.00}{2.00}{327}{110.00}{1.00}{328}
\put(146.00,81.00){\makebox(0,0)[cc]{$_{1^a}$}}
\put(175.00,81.00){\makebox(0,0)[cc]{$_{3^b}$}}
\put(206.00,81.00){\makebox(0,0)[cc]{$_{3^a}$}}
\put(115.00,81.00){\makebox(0,0)[cc]{$_{1^b}$}}
\put(176.00,65.00){\makebox(0,0)[cc]{$_{7^a}$}}
\put(205.00,65.00){\makebox(0,0)[cc]{$_{5^a}$}}
\put(145.00,65.00){\makebox(0,0)[cc]{$_{7^f}$}}
\put(116.00,65.00){\makebox(0,0)[cc]{$_{5^f}$}}
\put(146.00,49.00){\makebox(0,0)[cc]{$_{3^d}$}}
\put(175.00,49.00){\makebox(0,0)[cc]{$_{1^f}$}}
\put(115.00,49.00){\makebox(0,0)[cc]{$_{3^f}$}}
\put(206.00,49.00){\makebox(0,0)[cc]{$_{1^d}$}}
\put(176.00,33.00){\makebox(0,0)[cc]{$_{5^c}$}}
\put(205.00,33.00){\makebox(0,0)[cc]{$_{7^c}$}}
\put(145.00,33.00){\makebox(0,0)[cc]{$_{5^d}$}}
\put(116.00,33.00){\makebox(0,0)[cc]{$_{7^d}$}}
\put(146.00,17.00){\makebox(0,0)[cc]{$_{1^e}$}}
\put(175.00,17.00){\makebox(0,0)[cc]{$_{3^c}$}}
\put(115.00,17.00){\makebox(0,0)[cc]{$_{1^c}$}}
\put(206.00,17.00){\makebox(0,0)[cc]{$_{3^e}$}}
\put(176.00,3.00){\makebox(0,0)[cc]{$_{7^e}$}}
\put(205.00,3.00){\makebox(0,0)[cc]{$_{5^e}$}}
\put(145.00,3.00){\makebox(0,0)[cc]{$_{7^b}$}}
\put(116.00,3.00){\makebox(0,0)[cc]{$_{5^b}$}}
\emline{111.00}{1.00}{329}{114.00}{1.00}{330}
\emline{151.00}{1.00}{331}{148.00}{1.00}{332}
\emline{169.00}{1.00}{333}{173.00}{1.00}{334}
\emline{213.00}{1.00}{335}{209.00}{1.00}{336}
\end{picture}

\caption{A 6-hole and a star subgraph of $G$}
\end{figure}
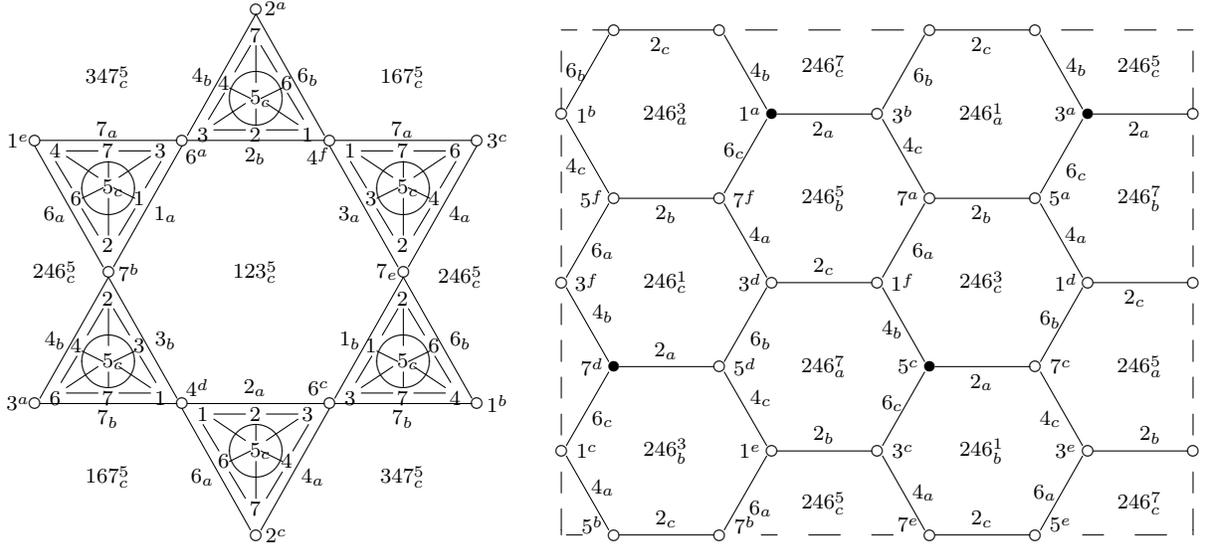

$G$ contains 84 6-holes obtainable from the octahedral triangles of $G$
based on these depictions of $\mathcal F$. This is exemplified on the
left side of Figure 4, where the upper-left triangle in Figure 3 appears as the bottom triangle, sharing its edge of weak color $2_a$ with the central 6-hole. The 6-cycle of weak colors associated to this 6-hole is $(2_a3_b1_a2_b3_a1_b)$.

Since the central color of each octahedral triangle having an edge in common with this 6-hole is $5_c$, and the line 123 has its points composing the weak colors of its edges, we denote this 6-hole by
$123_c^5$, as indicated in the figure. Similar denominations are given to the 6-holes neighboring this $123_c^5$ in the figure. We observe that the 6-holes neighboring $123_c^5$ on left and right coincide with $246_c^5$; on upper-left and lower-right with $347_c^5$; on lower-left and upper-right with $167_c^5$. These 4 6-holes form, together with the 8 octahedral triangles that appear in their generation (i.e. the 6 in the figure plus $(1^b,2^a,3^a)$ and $(1^e,2^c,3^c)$), a non-induced toroidal subgraph of $G$ that we may denote $[5]_c$.
In the same way, a non-induced subgraph $[w]_d$ is obtained, for each $w\in\mathcal F$ and $d\in\{a,b,c\}$. The subgraph $[[w]]_d$ induced by $[w]_d$ in $G$ is formed by its union with
the copies of $K_4$ in $G$ of the form $\langle x_ax_bx_c\rangle$, with $x_d=w$, a total of 6 copies of $K_4$ sharing each a 4-cycle with $[w]_d$. Each $[x_1x_bx_c]$ here is the union of such a 4-cycle plus two additional edges of weak color $w_d$. We get 21 subgraphs $[[w]]_d$ of $G$.

The 6-holes of the form $xyz_d^w$, where $xyz$ is a fixed line of $\mathcal F$, $d$ varies in $\{a,b,c\}$ and $w$ in $\mathcal F$, are 12 in number and conform a subgraph $[xyz]$ of $G$ isomorphic to the star Cayley graph $ST_4$, that can be defined as the graph with vertex set $S_4$ and each vertex $(a_0,a_1,a_2,a_3)\in S_4$ adjacent solely to $(a_1,a_0,a_2,a_3)$, $(a_2,a_1,a_0,a_3)$ and $(a_3,a_1,a_2,a_0)$; see \cite{akers,DS}. For example, the right side of Figure 4 depicts a (dotted) fundamental polygon of the torus whose convex hull contains a representation of the subgraph $[246]$ of $G$. This subgraphs $[xyz]$ are not induced  in $G$. However,
the graph $[[xyz]]$ induced in $G$ by each $[xyz]$
is the edge-disjoint union of $[xyz]$ with the edge-disjoint union of six
copies of $K_4$ in $G$, namely: $\langle xyz\rangle,\langle xzy\rangle ,\langle yxz\rangle,\langle yzx\rangle,\langle zxy\rangle,\langle zyx\rangle$.
Figure 4 has 4 vertices painted black, which span a copy of $K_4$ in $G$ but not in $[246]$.
We get 7 subgraphs $[xyz]$ of $G$.

Now, two new flag-transitive configurations associated to $G$ (apart from those cited in Section 6) are given by:
{\bf(a)} the 42 tetrahedra and 21 subgraphs $[w_d]$ in $G$, and
inclusion of a tetrahedron in a copy of $\overline{T}$; this is a
$(42_{3},21_{6})$-configuration;
{\bf(b)}
the 21 subgraphs $[w_d]$ and 7 copies of $ST_4$ in
$G$, and their sharing of a 6-hole; this is a
$(21_4,7_{12})$-configuration.

\newpage

When considered immersed in 3-space,
the 21 octahedra and 28 toroidal subgraphs of $G$ presented above have faces that appear in canceling pairs, allowing the visualization of a closed piecewise-linear 3-manifold.
We ask: which are the properties of this manifold?

\section{Open problems}

It remains to see whether $G$ is a Cayley graph or not.
On the other hand, the definition of $G$ may be extended by means
of projective planes, like the Fano plane, but
over larger fields than $GF(2)$, starting with $GF(3)$.
Moreover, the two conditions of the definition of $G$ in Section 3
may be taken to 3 conditions,
replacing ${\mathcal F}$ by a binary projective space $P(r-1,2)$ of dimension $r-1$, and
the Fano lines by subspaces of dimension $\sigma<r-1$,
where $2<r\in\ZZ$ and $\sigma\in(0,r-1)\cap\ZZ$,
and requiring, as a third condition, that
the points of intersection of a modified condition (b) form a projective hyperplane in
$P(r-1,2)$, (which was not required for $G$, since it was a ready conclusion).
The resulting graph, that appears in place of $G$, may not be even connected,
but the study of the component containing the lexicographically smallest vertex could
still be interesting. Another step would be taking the study over other fields, starting with the ternary one.

\vspace*{3mm}

\noindent{\bf Acknowledgement:} The author is grateful to Josep Rif\`a and Jaume Pujol,
for their friendship and
support, and to the referee, for his helpful observations, among which was suggesting the hemi-rhombicuboctahedron as the open neighborhood of each vertex of $G$.


\begin{thebibliography}{99}
\bibitem{akers} S. B. Akers and B. Krishnamurthy, {\it A group-theoretic model for symmetric interconnection networks}, IEEE Trans. Comput., 38(1989) 555--565.
\bibitem{biggs} N. Biggs, Algebraic Graph Theory, Cambridge University Press, 1993.
\bibitem{Cox} H.S.M. Coxeter, {\it Self-dual configurations and regular graphs},
Bull. Amer. Math. Soc., 56(1950) 413--455.
\bibitem{DS} I. J. Dejter and O. Serra, {\it Efficient dominating sets in Cayley graphs}, Discrete Applied Mathematics, 119 (2003) 319--328.
\bibitem{F} J. Folkman, {\it Regular line-symmetric graphs}, J. Combin. Theory, 3(1967), 215--232.
\bibitem{Gard} A. Gardiner, {\it Homogeneous graphs}, J. Combinatorial Theory (B), {\bf 20} (1976), 94-102.
\bibitem{I} D. C. Isaksen, C. Jankowski and S. Proctor, {\it On $K_*$-ultrahomogeneous graphs},
Ars Combinatoria, Volume LXXXII, (2007), 83--96.
\bibitem{Ronse} C. Ronse, {\it On homogeneous graphs}, J. London Math. Soc. (2) {\bf 17} (1978), 375--379.
\bibitem{Sheh} J. Sheehan, {\it Smoothly embeddable subgraphs}, J. London Math. Soc. (2)
{\bf 9} (1974), 212--218.
\end{thebibliography}
\end{document}